\documentclass[oneside,11pt]{article}

\usepackage[USenglish]{babel}
\usepackage[T1]{fontenc}
\usepackage{amsmath}
\usepackage{dsfont}
\usepackage{esint}
\usepackage{amsfonts}
\usepackage{amssymb}
\usepackage{amsthm}
\usepackage{mathrsfs}
\usepackage{graphicx}
\DeclareGraphicsExtensions{.eps}
\usepackage{cases}
\usepackage{amsmath}
\usepackage{units}
\usepackage{float}
\usepackage{braket}
\usepackage{subfigure}
\usepackage{rotating}
\usepackage{lineno}

\usepackage{xcolor}

\usepackage[small,justification=justified]{caption}

\newtheorem{remark}{Remark}

\newcommand{\R}{\mathbb{R}}
\newcommand{\dd}{\mathrm{d}}

\newcommand{\eps}{\varepsilon}

\usepackage[colorlinks=true, allcolors=blue]{hyperref}

\begin{document}

%%%%%%%%%%%%%%%%%%%%%%%%%%%%%%%%%%%%%%%%%%%%

\author{
Marzia Bisi\thanks{University of Parma, Italy.},
Andrea Bondesan\thanks{Gran Sasso Science Institute, Italy.}, %Corresponding author: \href{mailto: andrea.bondesan@gssi.it}{andrea.bondesan@gssi.it}},
Maria Groppi\footnotemark[1], 
Cinzia Soresina\thanks{University of Trento, Italy.}}

\title{Mesoscopic derivation of predator--prey dynamics via diffusive fast-reaction asymptotics: cross-diffusion and fear effects}

\maketitle

\date

\begin{abstract} \noindent
We propose a kinetic derivation of a spatial predator--prey model with cross-diffusion interactions and a Beddington--DeAngelis functional response. At the mesoscopic level, we divide the two populations into distinct behavioral states: prey can be active or hidden in a refuge, while predators can be either searching or handling. The microscopic state of each subpopulation is described in terms of group size and a continuous directional velocity variable, with state transitions governed by suitable Boltzmann-type interaction operators. By combining a diffusive scaling with a slow--fast formulation of the kinetic reaction operators, we formally derive a macroscopic reaction-cross-diffusion system for the total population densities of prey and predators. The limiting model consistently recovers the Beddington--DeAngelis functional response along with its associated nonlinear diffusion terms. Crucially, the kinetic upscaling yields a predator-dependent prey growth term taking a logistic form with an effective carrying capacity. This carrying capacity depends non-monotonically on predator density because prey redistribute between active and refuge states. Since density-dependent competition acts separately within these behavioral states, intermediate predator levels reduce aggregate crowding and thereby increase the effective carrying capacity of the total prey population. Overall, our derivation provides a mechanistic link connecting predator-induced behavioral shifts, refuge use, and macroscopic cross-diffusion.


\end{abstract}

%%%%%%%%%%%%%%%%%%%%%%%%%%%%%%%%%%%%%%%
%%%%    SECTION 1: INTRODUCTION    %%%%
%%%%%%%%%%%%%%%%%%%%%%%%%%%%%%%%%%%%%%%

\section{Introduction}

% Cross-diffusion
Spatial movement and behavioral responses to interspecific interactions are closely connected in ecological systems. In predator--prey communities, individuals may change their activity or habitat use in response to predation risk, prey availability, or feeding requirements. At the population level, these responses can produce spatial fluxes that depend on the densities of both species. Cross-diffusion systems provide a mathematical framework for describing such effects, allowing the movement of one population to depend on the spatial distribution of another. A central modeling question is how the resulting nonlinear diffusion terms can be related to specific ecological processes, and how these processes also affect the reaction terms governing population growth and trophic interactions.

% fast-reaction
A mechanism for deriving cross-diffusion consists in distinguishing behavioral states within each population and assuming that transitions between these states occur on a faster time scale than demographic processes~~\cite{IidaMimuraNinomiya2006, ConDesSor, DesSor, BrocchieriEtAl2021, EliasEtAl2022, iida2023cross, ConsoloEtAl2025, DesvillettesFiorentinoMautone2025, GianninoIuorioSoresina2026, brocchieri2026cross}. When individuals in different states have different movement characteristics, eliminating the fast processes leads to effective diffusion terms that depend on the local population densities. The same reduction also determines the effective reaction terms, since reproduction, competition, and predation act on populations whose distribution among behavioral states is constrained by the fast dynamics. This approach has been applied to predator--prey systems involving Holling-type~II and Beddington--DeAngelis functional responses~\cite{ConDesSor, DesSor, soresina2024fast}. In particular, following the mechanistic derivation of Geritz and Gyllenberg~\cite{GerGyl}, Desvillettes and Soresina~\cite{DesSor} obtained a non-triangular cross-diffusion system from a four-component reaction--diffusion model accounting for active and hidden prey, together with searching and handling predators. These derivations establish a connection between behavioral transitions, functional responses, and nonlinear spatial movement.

% kinetic derivation
Kinetic theory provides a complementary framework for connecting elementary interaction mechanisms with collective behaviors~\cite{bellomo2026scientific,BelHaOutYoo,FabLiaOut,GuaVerTos} in population dynamics. In this setting, a population is described through a distribution function over suitable microscopic variables, and interactions are encoded by operators (typically, Boltzmann and Fokker--Planck-type operators) specifying how these variables change. Macroscopic quantities are then recovered as moments of the distribution functions \cite{Per_book}. There is a vast literature of kinetic models for living systems in which the microscopic traits of individuals are characterized by their internal state, or activity \cite{BelBurLia,BurKno}, also in the context of population dynamics \cite{Lac}. For predator--prey interactions, Toscani and Zanella~\cite{TosZan} recently introduced an alternative description in which the microscopic variable represents the size of the interacting groups and the first moments of the mesoscopic solutions satisfy the Lotka--Volterra equations. Additional developments of this approach include the kinetic models studied in~\cite{BonMenTosZan, TosZan_eq}. Going further, spatial kinetic descriptions allow movement and interactions to be considered within the same framework, with effective diffusion equations obtained through suitable asymptotic scalings. We mention, for example, the recent work \cite{But} (whose approach is similar to ours, although based on the use of individuals' microscopic velocities for the kinetic variable, rather than the size of interacting groups), where predation is modeled at the kinetic level by a sequence of behavioral transitions and the use of slow--fast dynamics allows for the derivation of predator--prey standard diffusion models with Holling-type II and type III functional responses. Similarly, in the framework of reacting gas mixtures, Martal\`o and Travaglini~\cite{MarTra} derived a reaction--cross-diffusion system from kinetic equations, illustrating how cross-diffusion can also emerge from microscopic interactions occurring on different time scales.

% in this paper
In this paper, we develop a kinetic description of a spatial predator--prey system with a Beddington--DeAngelis functional response. Each species is divided into two behavioral states: prey can be active or hidden in a refuge, while predators can be searching or handling. Building on the group-size approach of~\cite{TosZan}, we describe each subpopulation by a distribution function and represent the relevant interactions through Boltzmann-type operators. The interaction rules account for prey capture, predator-induced refuge use, return to activity, transitions between searching and handling, and competition within each prey state. Birth and death processes are included through additional linear drift operators. We then introduce spatial transport and relaxation operators to describe movement within each behavioral state.

By combining a diffusive scaling with a slow--fast formulation of the kinetic reaction operators, we formally derive a macroscopic reaction--cross-diffusion system for the total prey and predator densities. The scaling distinguishes dominant relaxation, fast behavioral transitions, and slow demographic processes, thereby making explicit the hierarchy of time scales underlying the reduction. In particular, convergence to local kinetic equilibrium for each distribution provides at leading order a link between the corresponding flux and concentration gradient, while the behavioral transitions impose at intermediate scales additional relations between the state densities, allowing going beyond the usual quadratic terms prescribed by Boltzmann operators and recovering more complex nonlinearities in the combined diffusion--fast-reaction limit. Specifically, the limiting system recovers the Beddington--DeAngelis functional response and the nonlinear diffusion structure obtained in~\cite{DesSor}, together with a novel prey-growth term determined by the demographic assumptions at the kinetic level. As a byproduct, if no state dichotomy in the prey is assumed, we are also able to recover the Holling-type II structure of \cite{ConDesSor}. The derivation thus relates the macroscopic coefficients and nonlinearities to the underlying interaction rules and clarifies how the behavioral mechanisms contribute to both spatial movement and population dynamics. To our knowledge, this is the first derivation of this kind for cross-diffusion predator--prey models. Moreover, in contrast to previous works on reacting gas mixtures \cite{MarTra}, where a nonstandard rescaling of time and space is introduced to bring the linear relaxation operator to order $\mathcal{O}(\eps^2)$, rather than the usual $\mathcal{O}(\eps)$ order in diffusive limits (thus leading to a Hilbert expansion with no first-order kinetic correction), we show that cross-diffusion emerges under a standard parabolic rescaling of space and time, allowing us to close the moment equations already at the first-order kinetic correction of the Hilbert expansion.

% sottoprodotto
An additional outcome concerns the effect of predator-induced refuge use on prey growth. When active and hidden prey have the same intrinsic growth rate, and competition acts separately within the two behavioral states, the limiting growth term takes a logistic form with a predator-dependent effective carrying capacity. This dependence is non-monotone: an intermediate predator density distributes prey more evenly between the active and hidden states, reducing aggregate competition and increasing the effective carrying capacity. At high predator densities, prey become concentrated in the hidden state, and this competitive release diminishes. The resulting carrying capacity therefore reflects the redistribution of prey and the structure of intraspecific competition.

% fear effect
This observation also provides a connection with the non-consumptive effects of predators. Predator-induced changes in prey behavior and habitat use are well documented~\cite{lima1990behavioral,werner1983experimental,werner2003trait}, and their reproductive costs have motivated mathematical models incorporating the so-called fear effect~\cite{wang2016modelling}. In our framework, assuming that hidden prey temporarily cease to reproduce yields an effective intrinsic growth rate that decreases with predator density. As far as we know, this work provides the first kinetic derivation of the reproductive cost associated with the fear effect. The predator-dependent reduction in prey reproduction emerges from fast transitions between active, reproducing prey and hidden, non-reproducing prey. The same behavioral redistribution also modifies the density-dependent component of prey growth, with consequences that depend on whether hidden individuals remain subject to competition. The kinetic derivation therefore provides a common mechanistic description of the reproductive cost of refuge use and its effects on intraspecific competition, showing how both contribute to the effective prey-growth term.

% Paper organisation
The paper is organized as follows. Section~\ref{sec:kinmod} introduces the interaction mechanisms and the kinetic model, and derives the corresponding equations for the population densities. Section~\ref{sec:diff} presents the spatial kinetic description and the formal diffusive fast-reaction limit leading to the macroscopic reaction--cross-diffusion system. Section~\ref{sec:ecological_interpretation} examines the predator-dependent prey-growth term and discusses its ecological interpretation, including its connection with fear effects. Finally, Section~\ref{sec:conclusions} summarizes the main findings and outlines directions for future research.

%%%%%%%%%%%%%%%%%%%%%%%%%%%%%%%%%%%%%%%%%%%%%%%%%%%%%
%%%%    SECTION 2: PHYSICAL SETTING AND MODEL    %%%%
%%%%%%%%%%%%%%%%%%%%%%%%%%%%%%%%%%%%%%%%%%%%%%%%%%%%%

\section{Physical setting and kinetic model}\label{sec:kinmod}
 
We consider a population of individuals composed of prey and predators, where each species is characterized by two different behavioral states. Prey are either active~($P$), engaging in foraging and mating, or hidden~($R$), taking refuge from predators. Predators are either searching~($S$), actively hunting for prey, or handling~($H$), a state that includes eating, digesting, and resting. The set of behavioral states is therefore given by~$\mathcal{S} = \{P,R,S,H\}$. Following the kinetic formalism introduced in \cite{BonMenTosZan, TosZan}, we describe the evolution of each behavioral state $J \in \mathcal{S}$ using a distribution function $f_J = f_J(t,z)$ that depends on time $t \in \R_+$ and on a kinetic variable $z \in \R_+$, characterizing the size of the considered group of individuals. In other words, the quantity $f_J(t,z)\,\dd z$ provides the number of groups of individuals belonging to state $J \in \mathcal{S}$ and having size in $[z, z+\dd z]$. We then associate with each distribution function $f_J$ its first-order macroscopic moment in~$z \in \R_+$
\begin{equation*}
    m_J(t) = \int_{\R_+} z f_J(t,z) \dd z,  
\end{equation*}
defining the number of individuals (or density) of state $J \in \mathcal{S}$. We also assume that both species can interact with a fixed surrounding environment of resources (additional food, shelters, etc.) described by a time-independent distribution function $f_E = f_E(z)$ satisfying the normalization
\begin{equation*}
    \int_{\R_+} z f_E(z) \dd z = 1.
\end{equation*}

\begin{remark}[Macroscopic observables] \label{rem:groups vs individuals 1}
    In our setting, it is essential to distinguish between the \textit{number of groups} and the \textit{number of individuals} in each behavioral state. Unlike kinetic descriptions in which the zeroth-order moment gives the population density, here the macroscopic quantities
    \begin{equation*}
        n_J(t) = \int_{\R_+} f_J(t,z) \dd z 
    \end{equation*} 
    represent the number of groups in each state $J \in \mathcal{S}$, whereas the first-order moments $m_J$ give the corresponding population densities. This interpretation follows directly from the meaning of the microscopic variable $z \in \R_+$, which denotes group size, and is reminiscent of the kinetic framework for the Smoluchowski equation~\cite{Smo}.
    As a consequence, the number of groups and the number of individuals in a certain state need not satisfy the same conservation laws. Births and deaths, for instance, may change group sizes while leaving the number of groups unchanged, as new or remaining individuals are distributed among existing groups. Conversely, groups may merge or split, changing their number without affecting the total number of individuals in a given state. This distinction will be crucial in the construction of the kinetic operators in the following.
\end{remark}

To formulate the kinetic system governing the evolution of the distribution functions $f_J$, $J \in \mathcal{S}$, we first describe the microscopic processes that modify group sizes and transfer individuals between behavioral states. We assume that both active prey $P$ and hidden prey $R$ may reproduce through access to environmental resources, at rates independent of predator presence. Searching predators $S$ and handling predators $H$ are subject to natural mortality due to aging, disease, competition, and other environmental factors, independently of predation. In addition to these demographic processes, individuals may change behavioral state through interactions. For predators, these transitions are governed by the following four microscopic processes:
\begin{align}
    & S + P \rightarrow H, \label{eq:S+P to H} \\[2mm]
    & S + P \rightarrow H + P, \label{eq:S+P to H+P} \\[2mm]
    & H + E \rightarrow H + S, \label{eq:H+E to H+S} \\[2mm]
    & H + E \rightarrow S. \label{eq:H+E to S}
\end{align}
Reactions~\eqref{eq:S+P to H} and~\eqref{eq:S+P to H+P} describe the transition from searching to handling predators following a predation attempt. A successful capture removes the prey from the system, whereas an unsuccessful attempt leads to a temporary resting or recovery phase, which is also included in the handling state $H$. Handling predators may reproduce through reaction~\eqref{eq:H+E to H+S} or return to the searching state through reaction~\eqref{eq:H+E to S}. Both processes are modeled as interactions with the surrounding environment. On the other hand, prey's behavioral states are involved in the five interactions
\begin{align}
    & P + S \rightarrow R + S, \label{eq:P+S to R+S} \\[2mm]
    & P + H \rightarrow R + H, \label{eq:P+H to R+H} \\[2mm]
    & R + E \rightarrow P, \label{eq:R+E to P} \\[2mm]
    & P + P \rightarrow P, \label{eq:P+P to P} \\[2mm]
    & R + R \rightarrow R. \label{eq:R+R to R}
\end{align}
Processes~\eqref{eq:P+S to R+S} and~\eqref{eq:P+H to R+H} describe active prey taking refuge after escaping from a searching or handling predator. Hidden prey may subsequently return to the active state through interactions with the surrounding environment, as described by process~\eqref{eq:R+E to P}. Finally, competition within the active and hidden prey states may lead to mortality, as modeled by processes~\eqref{eq:P+P to P} and~\eqref{eq:R+R to R}, respectively.

%%%%    SUBSECTION: KINETIC EQUATIONS    %%%%
\subsection{The kinetic system}

For notational convenience, let us introduce the vector distribution function $\mathbf{f} = (f_P,f_R,f_S,f_H,f_E)$, where we recall that $f_E(z)$ does not change over time. Based on the previous discussion, the kinetic system governing the evolution of $\mathbf{f}$ reads
\begin{equation} \label{eq:kinetic model}
    \frac{\partial}{\partial t} f_J = \mathcal{R}_J(f_J) + \mathcal{Q}_J(\mathbf{f},\mathbf{f}), \quad t \in \R_+, \, z \in \R_+, \, J \in \mathcal{S},
\end{equation}
where the redistribution-like operators $\mathcal{R}_J(f_J)$ characterize the natural birth and death processes to which prey and predators are respectively subjected, while the Boltzmann-like operators $\mathcal{Q}_J(\mathbf{f},\mathbf{f})$ account for all the interactions \eqref{eq:S+P to H}--\eqref{eq:R+R to R} prescribing the switches between different behavioral states.

More precisely, the linear operators $\mathcal{R}_J(f_J)$, $J \in \mathcal{S}$, are inspired by similar ones introduced in \cite{BisSpiTos} to describe the redistribution of wealth collected by taxation in a simple market economy. In particular, the birth operators for prey $P$ and $R$ take the form of a drift
\begin{equation} \label{eq:birth operators}
    \mathcal{R}_J(f_J)(t,z) = -\alpha_J \frac{\partial}{\partial z} \big( z f_J(t,z) \big), \quad J = P,R,
\end{equation}
where the coefficients $\alpha_P, \alpha_R > 0$ measure the baseline birth rates of active and hidden prey. Notice that in general it is reasonable to assume different birth rates for active and hidden prey, as the first ones are more likely to reproduce thanks to energy trade-offs (active individuals spend more time foraging and consuming high-quality food, giving them more energy to allocate toward producing offspring, while hidden individuals in refuge conserve energy but have highly restricted access to food, leading to lower nutritional status and a lower birth rate), mating opportunities (active individuals move around the habitat, greatly increasing their chances of encountering, courting, and successfully mating with partners, while individuals hidden away in refuges have severely limited social or reproductive interactions), and stress-induced infertility (staying in a long-term hiding state is often a response to intense predator pressure, and high stress levels elevate hormones that can suppress the reproductive systems of hidden individuals, lowering their birth rate or causing them to skip breeding cycles entirely). The death operators for predators $S$ and $H$ have the same form but are taken with opposite signs, namely
\begin{equation} \label{eq:death operators}
    \mathcal{R}_J(f_J)(t,z) = \mu \frac{\partial}{\partial z} \big( z f_J(t,z) \big), \quad J = S, H,
\end{equation}
and now the coefficient $\mu > 0$ defines the baseline death rate of predators. Taking the same death rate for both searching and handling predators is justified, for example, by assuming that the handling phase is strictly metabolic (hence, it happens on a short time scale) and that predators are not menaced by other hunting species. 

On the other hand, the bilinear operators $\mathcal{Q}_J(\mathbf{f},\mathbf{f})$, $J \in \mathcal{S}$, characterizing the reactions \eqref{eq:S+P to H}--\eqref{eq:R+R to R} split into several components that we write here explicitly for easier visualization. Active prey $P$ are involved in the five processes \eqref{eq:S+P to H}, \eqref{eq:P+S to R+S}, \eqref{eq:P+H to R+H}, \eqref{eq:R+E to P}, and \eqref{eq:P+P to P}, hence the corresponding Boltzmann operator $\mathcal{Q}_P(\mathbf{f},\mathbf{f})$ is given by
\begin{equation} \label{eq:Q_P}
    \begin{split}
        \mathcal{Q}_P(\mathbf{f},\mathbf{f}) &= \mathcal{Q}_P^{S + P \rightarrow H}(f_S,f_P) + \mathcal{Q}_P^{P + S \rightarrow R + S}(f_P,f_S) + \mathcal{Q}_P^{P + H \rightarrow R + H}(f_P,f_H) \\[4mm]
        &\qquad + \mathcal{Q}_P^{R + E \rightarrow P}(f_R,f_E) + \mathcal{Q}_P^{P + P \rightarrow P}(f_P,f_P).
    \end{split}
\end{equation}
The number of prey $R$ who found a refuge is modified by the four processes \eqref{eq:P+S to R+S},~\eqref{eq:P+H to R+H},~\eqref{eq:R+E to P}, and~\eqref{eq:R+R to R}, therefore~$\mathcal{Q}_R(\mathbf{f},\mathbf{f})$ writes as
\begin{equation} \label{eq:Q_R}
    \begin{split}
        \mathcal{Q}_R(\mathbf{f},\mathbf{f}) &= \mathcal{Q}_R^{P + S \rightarrow R + S}(f_P,f_S) + \mathcal{Q}_R^{P + H \rightarrow R + H}(f_P,f_H) + \mathcal{Q}_R^{R + E \rightarrow P}(f_R,f_E) + \mathcal{Q}_R^{R + R \rightarrow R}(f_R,f_R).
    \end{split}
\end{equation}
Searching predators $S$ are subject to the four processes \eqref{eq:S+P to H}, \eqref{eq:S+P to H+P}, \eqref{eq:H+E to H+S}, and~\eqref{eq:H+E to S}, leading to a Boltzmann operator~$\mathcal{Q}_S(\mathbf{f},\mathbf{f})$ of the form
\begin{equation} \label{eq:Q_S}
    \begin{split}
        \mathcal{Q}_S(\mathbf{f},\mathbf{f}) &= \mathcal{Q}_S^{S + P \rightarrow H}(f_S,f_P) + \mathcal{Q}_S^{S + P \rightarrow H + P}(f_S,f_P) + \mathcal{Q}_S^{H + E \rightarrow H + S}(f_H,f_E) + \mathcal{Q}_S^{H + E \rightarrow S}(f_H,f_E).
    \end{split}
\end{equation}
Finally, changes in the number of handling predators $H$ are affected by the three processes \eqref{eq:S+P to H}, \eqref{eq:S+P to H+P}, and~\eqref{eq:H+E to S}, and in this case the Boltzmann operator $\mathcal{Q}_H(\mathbf{f},\mathbf{f})$ reads
\begin{equation} \label{eq:Q_H}
    \begin{split}
        \mathcal{Q}_H(\mathbf{f},\mathbf{f}) &= \mathcal{Q}_H^{S + P \rightarrow H}(f_S,f_P) + \mathcal{Q}_H^{S + P \rightarrow H + P}(f_S,f_P) + \mathcal{Q}_H^{H + E \rightarrow S}(f_H,f_E).
    \end{split}
\end{equation}

In order to write the operators \eqref{eq:Q_P}--\eqref{eq:Q_H}, we start by presenting an abstract setting that will be useful to introduce the key modeling ideas behind their formulation. We will then exploit these considerations to describe the explicit structure of the different components of $\mathcal{Q}_J(\mathbf{f},\mathbf{f})$, for each process \eqref{eq:S+P to H}--\eqref{eq:R+R to R} separately.

\medskip
\noindent \textbf{General structure of the bilinear operators.} Let us begin by considering an abstract transition between ordered states $I + K \rightarrow I' + K'$, with $I \neq I'$ and $K \neq K'$. The interpretation we give to this generic transition is the following: when two groups in states $I$ and $K$ interact, a fraction of their individuals switch state to~$I'$ and~$K'$ respectively, and the number of prey and predators from the pre-interaction states is redistributed among groups. Then, the Boltzmann operators $\mathcal{Q}_J^{I + K \rightarrow I' + K'}(f_I,f_K)$, $J \in \{I, I', K, K'\}$, associated with the process $I + K \rightarrow I' + K'$, act only on the size variable $z \in \R_+$ and are thus naturally defined in weak form by testing them against regular functions $\varphi(z)$. Their explicit structure reads
\begin{equation} \label{eq:collision operators}
    \begin{split}
        \int_{\R_+} \varphi(z) \mathcal{Q}_J^{I + K \rightarrow I' + K'}&(f_I,f_K)(z) \dd z = \\[2mm]
        & \int_{\R_+ \times \R_+} \kappa_J^{I + K \rightarrow I' + K'}(z,z_*) \mathbb{E}\big(\varphi(z')-\delta_{JI} \varphi(z)\big) f_I(z) f_K(z_*) \dd z \dd z_*, \quad J = I, I', \\[6mm]
        & \int_{\R_+ \times \R_+} \kappa_J^{I + K \rightarrow I' + K'}(z,z_*) \mathbb{E}\big(\varphi(z')-\delta_{JK} \varphi(z)\big) f_I(z_*) f_K(z) \dd z \dd z_*, \quad J = K, K',
    \end{split}
\end{equation}
corresponding to a loss and redistribution of individuals from states $I$ and $K$, and a gain of individuals from states $I'$ and $K'$, respectively. In \eqref{eq:collision operators}, the functions $\kappa_J^{I + K \rightarrow I' + K'}(z,z_*) \geq 0$ denote the transition kernels of individuals from states $(I,K)$ to states $(I',K')$, while $\delta_{JJ'} = \mathds{1}_{\{J = J'\}}$ denotes the Kronecker delta function. Moreover, the post-interaction sizes $z' \in \R_+$ appearing in the integrands are determined by the microscopic updates
\begin{equation} \label{eq:microsopic interactions}
    \begin{array}{lll}
        z' = z - \psi_K^{I + K \rightarrow I' + K'}(z_*) z + \eta_K^{I + K \rightarrow I' + K'}(z_*) z, & \ \text{if} & J = I, \\[6mm]
        z' = z - \psi_I^{I + K \rightarrow I' + K'}(z_*) z + \eta_I^{I + K \rightarrow I' + K'}(z_*) z, & \ \text{if} & J = K, \\[6mm]
        z' = \psi_K^{I + K \rightarrow I' + K'}(z_*) z + \tilde{\eta}_K^{I + K \rightarrow I' + K'}(z_*) z, & \ \text{if} & J = I', \\[6mm]
        z' = \psi_I^{I + K \rightarrow I' + K'}(z_*) z + \tilde{\eta}_I^{I + K \rightarrow I' + K'}(z_*) z, & \ \text{if} & J = K'.
    \end{array}
\end{equation}
The microscopic interaction rules \eqref{eq:microsopic interactions} describe a redistribution mechanism whereby the size of a group initially belonging to state $I$ (respectively, $K$) is divided between its original state and the corresponding destination state $I'$ (respectively, $K'$). The (deterministic) state-transition efficiency functions $0 \leq \psi_K^{I + K \rightarrow I' + K'} < 1$ and $0 \leq \psi_J^{I + K \rightarrow I' + K'} < 1$ quantify the fractions of size transferred to the destination states, while the perturbations $\eta_J^{I + K \rightarrow I' + K'}(z_*)$ and $\tilde{\eta}_J^{I + K \rightarrow I' + K'}(z_*)$ are zero-mean random variables accounting for unpredictable fluctuations in this redistribution. In particular, the zero-mean assumptions on the random variables ensure that the total size is conserved in expectation, although conservation is not necessarily satisfied for each individual realization. To guarantee the nonnegativity of the post-interaction sizes $z'$ in the transformations of \eqref{eq:microsopic interactions}, we assume that, for $J = I,K$,
\begin{equation*}
    \begin{split}
        \eta_J^{I + K \rightarrow I' + K'}(z_*) &\geq -\left(1 - \psi_J^{I + K \rightarrow I' + K'}(z_*)\right), \\[4mm]
        \tilde{\eta}_J^{I + K \rightarrow I' + K'}(z_*) &\geq - \psi_J^{I + K \rightarrow I' + K'}(z_*),
    \end{split}
\end{equation*}
almost surely. Moreover, note that the expected value appearing in the integral \eqref{eq:collision operators} is taken with respect to the probability density functions associated with $\eta_J^{I + K \rightarrow I' + K'}(z_*)$ and $\tilde{\eta}_J^{I + K \rightarrow I' + K'}(z_*)$.

\begin{remark}[Example of deterministic interaction functions] \label{rem:positivity of post-interactions}
    A relevant modeling example for the functions~$\psi_J^{I + K \rightarrow I' + K'}(z_*)$ is given by
    \begin{equation*}
        \psi_J^{I + K \rightarrow I' + K'}(z_*) = \psi_{IK}^{I'K'} \frac{z_*}{1+z_*},
    \end{equation*}
    where the constant $\psi_{IK}^{I'K'} \in (0,1)$ defines the baseline interaction intensity of process $I + K \rightarrow I' + K'$ (see~\cite{BonMenTosZan} for more details). Specifically, this Holling-type II interaction describes a saturation phenomenon in the effectiveness of the state-transition mechanism, with $\psi_{IK}^{I'K'}$ representing its maximum efficiency.
\end{remark}

The four post-interaction sizes $z'$ in \eqref{eq:microsopic interactions} are identified by the type of process $I + K \rightarrow I' + K'$ and by the behavioral state $J \in \mathcal{S}$ involved in the process. The first and second transformations prescribe a loss of individuals from the pre-interaction states $I$ and $K$, while the third and fourth ones cover the complementary event of individuals gained by the post-interaction states $I'$ and $K'$. Considering, for example, the case~$J = I$, prey or predators of size $z$ belonging to state $I$ lose a fraction of their group $-\psi_K^{I + K \rightarrow I' + K'}(z_*)z$ proportional to the number of individuals from state $K$ that are involved in the microscopic interaction, and the remaining ones $z'$ are simultaneously redistributed within the corresponding species. Conversely, when~$J = I'$ in the third transformation, prey or predators from state $I'$ gain a new group of size $z'$ given by the fraction~$\psi_K^{I + K \rightarrow I' + K'}(z_*)z$ of individuals lost from state $I$ in the first update. Cases~$J = K, \,K'$ work identically.

\begin{remark}[Follow-up on Remark \ref{rem:groups vs individuals 1}] \label{rem:groups vs individuals 2}
    The presence of the Kronecker delta inside relations \eqref{eq:collision operators} allows us to distinguish between a full gain/loss structure ($\delta_{JJ'} = 1$) of the Boltzmann operators for the pre-interaction states $I, K$ and a pure gain term ($\delta_{JJ'} = 0$) for the post-interaction states $I',\, K'$. The first case encodes at the same time a loss and a redistribution of individuals, while the latter leads to an overall increase in the number of groups. This fact stems directly from the conservation properties of the operators $Q_J^{I + K \rightarrow I' + K'}(f_I,f_K)$, since choosing $\varphi(z)=1$ in \eqref{eq:collision operators} one can see that
    \begin{equation*}
        \int_{\R_+} \mathcal{Q}_J^{I + K \rightarrow I' + K'}(f_I,f_K)(z) \dd z = 
        \left\{ \begin{array}{ll}
            0, & J = I,K, \\[4mm]
            \displaystyle \int_{\R_+ \times \R_+} \kappa_J^{I + K \rightarrow I' + K'}(z,z_*) f_I(z) f_K(z_*) \dd z \dd z_* \geq 0, & J = I', \\[6mm]
            \displaystyle \int_{\R_+ \times \R_+} \kappa_J^{I + K \rightarrow I' + K'}(z,z_*) f_I(z_*) f_K(z) \dd z \dd z_* \geq 0, & J = K'.
        \end{array} \right.
    \end{equation*}
    In particular, although states $I$ and $K$ have lost a fraction of their individuals, the number of groups $n_I$ and~$n_K$ are preserved during the interaction, meaning that individuals have redistributed to form new groups. On the contrary, the process $I + K \rightarrow I' + K'$ produces a net gain of both individuals and groups in states $I'$ and~$K'$, since in this case the macroscopic quantities $n_{I'}$ and $n_{K'}$ tend to increase. Additionally, we point out that our choice of modeling the action of the Boltzmann operators differently for the pre- and post-interaction states (gain/loss and redistribution for $J = I,\, K$ and pure gain for $J = I',\, K'$) is crucial from a mathematical standpoint, to ensure that the flow of system \eqref{eq:kinetic model} preserves the nonnegativity of the distributions $f_J$.
\end{remark}

Using the ideas behind the construction of the collision operators \eqref{eq:collision operators}, we are now ready to present the explicit structure of the Boltzmann operators \eqref{eq:Q_P}--\eqref{eq:Q_H}, by studying processes \eqref{eq:S+P to H}--\eqref{eq:R+R to R} one by one.

\medskip
\noindent \textbf{Process $\mathbf{S + P \rightarrow H}$.}
Process \eqref{eq:S+P to H} exemplifies a situation in which one of the pre-interaction states is transformed while the other is lost. Specifically, when a searching predator $S$ successfully catches an active prey $P$, the latter dies, and the predator can rest, becoming handling $H$ for some time. Therefore, this process involves three different Boltzmann operators accounting for the losses of prey $P$ and predators $S$, and for the gain and redistribution of predators $H$. These are given by
\begin{equation*}
    \begin{split}
        \int_{\R_+} \varphi(z) Q_J^{S + P \rightarrow H}&(f_S,f_P)(z) \dd z = \\[2mm]
        & \int_{\R_+ \times \R_+} \kappa_S^{S + P \rightarrow H}(z,z_*) \mathbb{E}\big(\varphi(z')-\varphi(z)\big) f_S(z) f_P(z_*) \dd z \dd z_*, \quad J = S, \\[3mm]
        & \int_{\R_+ \times \R_+} \kappa_P^{S + P \rightarrow H}(z,z_*) \mathbb{E}\big(\varphi(z')-\varphi(z)\big) f_S(z_*) f_P(z) \dd z \dd z_*, \quad J = P, \\[3mm]
        & \int_{\R_+ \times \R_+} \kappa_H^{S + P \rightarrow H}(z,z_*) \mathbb{E}\big(\varphi(z')\big) f_S(z) f_P(z_*) \dd z \dd z_*, \quad J = H,
    \end{split}
\end{equation*}
and the post-interaction size $z'$ is determined for each of the three behavioral states involved as
\begin{equation*}
    \begin{array}{ll}
        z' = z - \psi_P^{S + P \rightarrow H}(z_*) z + \eta_P^{S + P \rightarrow H}(z_*) z, & J = S, \\[4mm]
        z' = z - \psi_S^{S + P \rightarrow H}(z_*) z + \eta_S^{S + P \rightarrow H}(z_*) z, & J = P, \\[4mm]
        z' = \psi_P^{S + P \rightarrow H}(z_*) z + \tilde{\eta}_P^{S + P \rightarrow H}(z_*) z, & J = H.
    \end{array}
\end{equation*}

\medskip
\noindent \textbf{Process $\mathbf{S + P \rightarrow H + P}$.} In process \eqref{eq:S+P to H+P}, one of the pre-interaction states is transformed, while the other is kept fixed. Specifically, when searching predators $S$ hunt for active prey $P$ but fail to catch them, they move into a resting state $H$ before going back to hunting for new prey. In particular, we assume that prey remain in the active state $P$ even after the interaction. Hence, the process involves two separate Boltzmann operators that account for a loss and redistribution of searching predators $S$, and a gain of handling ones $H$. These read
\begin{equation*}
    \begin{split}
        \int_{\R_+} \varphi(z) Q_J^{S + P \rightarrow H + P}&(f_S,f_P)(z) \dd z = \\[2mm]
        & \int_{\R_+ \times \R_+} \kappa_S^{S + P \rightarrow H + P}(z,z_*) \mathbb{E}\big(\varphi(z')-\varphi(z)\big) f_S(z) f_P(z_*) \dd z \dd z_*, \quad J = S, \\[4mm]
        & \int_{\R_+ \times \R_+} \kappa_H^{S + P \rightarrow H + P}(z,z_*) \mathbb{E}\big(\varphi(z')\big) f_S(z) f_P(z_*) \dd z \dd z_*, \quad J = H,
    \end{split}
\end{equation*}
and the post-interaction sizes $z'$ are given by
\begin{equation*}
    \begin{array}{ll}
        z' = z - \psi_P^{S + P \rightarrow H + P}(z_*) z + \eta_P^{S + P \rightarrow H + P}(z_*) z, & J = S, \\[4mm]
        z' = \psi_P^{S + P \rightarrow H + P}(z_*) z + \tilde{\eta}_P^{S + P \rightarrow H + P}(z_*) z, & J = H.
    \end{array}
\end{equation*}

\medskip
\noindent \textbf{Process $\mathbf{H + E \rightarrow H + S}$.} The third process \eqref{eq:H+E to H+S} involves instead an interaction of handling predators $H$ with the fixed environment $E$ (which is not modified in the process), allowing them to reproduce $H$ and go back to a searching state $S$. We thus implicitly assume that the newborn group of handling predators matches the pre-interaction one, to ensure that only searching predators are produced in the process. In particular, we deduce that $Q_H^{H + E \rightarrow H + S}(f_H,f_E) = 0$, and the only nonzero collision operator is the one accounting for a gain of searching predators $S$ given by
\begin{equation*}
    \begin{split}
        \int_{\R_+} \varphi(z) Q_S^{H + E \rightarrow H + S}(f_H,f_E)(z) \dd z = \int_{\R_+ \times \R_+} \kappa_S^{H + E \rightarrow H + S}(z,z_*) \mathbb{E}\big(\varphi(z')\big) f_H(z) f_E(z_*) \dd z \dd z_*,
    \end{split}
\end{equation*}
where the post-interaction size $z'$ reads
\begin{equation*}
    \begin{split}
        z' = \psi_E^{H + E \rightarrow H + S}(z_*) z + \tilde{\eta}_E^{H + E \rightarrow H + S}(z_*) z.
    \end{split}
\end{equation*}
Notice that the operator is linear in this case, since the distribution of resources $f_E(z_*)$ does not evolve. Obviously, the same happens for all processes involving an interaction with the background environment.

\medskip
\noindent \textbf{Process $\mathbf{H + E \rightarrow S}$.}
Process \eqref{eq:H+E to S} describes the transition of predators from a handling state $H$ back to a searching one $S$, as the result of an interaction with the fixed environment $E$. The construction of the Boltzmann operators in this case works in the same way as for process $S + P \rightarrow H$, except that here there is no actual loss of individuals from state $E$ and the only two resulting operators for the balance of states $H$ and $S$ are linear. Their structure explicitly reads
\begin{equation*}
    \begin{split}
        \int_{\R_+} \varphi(z) Q_J^{H + E \rightarrow S}&(f_H,f_E)(z) \dd z = \\[2mm]
        & \int_{\R_+ \times \R_+} \kappa_H^{H + E \rightarrow S}(z,z_*) \mathbb{E}\big(\varphi(z')-\varphi(z)\big) f_H(z) f_E(z_*) \dd z \dd z_*, \quad J = H, \\[4mm]
        & \int_{\R_+ \times \R_+} \kappa_S^{H + E \rightarrow S}(z,z_*) \mathbb{E}\big(\varphi(z')\big) f_H(z) f_E(z_*) \dd z \dd z_*, \quad J = S,
    \end{split}
\end{equation*}
where the post-interaction sizes $z'$ are given by
\begin{equation*}
    \begin{array}{ll}
        z' = z - \psi_E^{H + E \rightarrow S}(z_*) z + \eta_E^{H + E \rightarrow S}(z_*) z, & J = H, \\[4mm]
        z' = \psi_E^{H + E \rightarrow S}(z_*) z + \tilde{\eta}_E^{H + E \rightarrow S}(z_*) z, & J = S.
    \end{array}
\end{equation*}

\medskip
\noindent \textbf{Processes $\mathbf{P + S \rightarrow R + S}$ and $\mathbf{P + H \rightarrow R + H}$.} To construct the collision operators modeling processes~\eqref{eq:P+S to R+S} and \eqref{eq:P+H to R+H}, we follow the same reasoning used for the transition $S+P \rightarrow H+P$. Here, when active prey~$P$ encounter either searching $S$ or handling $H$ predators, they seek refuge $R$ to hide from them and lower their risk of being eaten. In particular, we assume that predators remain in the hunting $S$ or resting $H$ state even after the interaction. Hence, both processes involve two separate Boltzmann operators that account for a loss and redistribution of prey $P$, and a gain of prey $R$. For process $P + S \rightarrow R + S$, these read
\begin{equation*}
    \begin{split}
        \int_{\R_+} \varphi(z) Q_J^{P + S \rightarrow R + S}&(f_P,f_S)(z) \dd z = \\[2mm]
        & \int_{\R_+ \times \R_+} \kappa_P^{P + S \rightarrow R + S}(z,z_*) \mathbb{E}\big(\varphi(z')-\varphi(z)\big) f_P(z) f_S(z_*) \dd z \dd z_*, \quad J = P, \\[4mm]
        & \int_{\R_+ \times \R_+} \kappa_R^{P + S \rightarrow R + S}(z,z_*) \mathbb{E}\big(\varphi(z')\big) f_P(z) f_S(z_*) \dd z \dd z_*, \quad J = R,
    \end{split}
\end{equation*}
with post-interaction sizes $z'$ given by
\begin{equation*}
    \begin{array}{ll}
        z' = z - \psi_S^{P + S \rightarrow R + S}(z_*) z + \eta_S^{P + S \rightarrow R + S}(z_*) z, & J = P, \\[4mm]
        z' = \psi_S^{P + S \rightarrow R + S}(z_*) z + \tilde{\eta}_S^{P + S \rightarrow R + S}(z_*) z, & J = R.
    \end{array}
\end{equation*}
The collision operators and post-interaction sizes for process $P + H \rightarrow R + H$ are written in the same way, by simply replacing each instance of $S$ with $H$.

\medskip
\noindent \textbf{Process $\mathbf{R + E \rightarrow P}$.}
Process \eqref{eq:R+E to P} provides the equivalent of transition $H+E \rightarrow S$ for prey populations, where hidden prey $R$ go back to the active state $P$ after interacting with the environment $E$. Therefore, the two Boltzmann operators accounting for the loss of prey $R$ and the gain of prey $P$ display a similar structure, namely
\begin{equation*}
    \begin{split}
        \int_{\R_+} \varphi(z) Q_J^{R + E \rightarrow P}&(f_R,f_E)(z) \dd z = \\[2mm]
        & \int_{\R_+ \times \R_+} \kappa_R^{R + E \rightarrow P}(z,z_*) \mathbb{E}\big(\varphi(z')-\varphi(z)\big) f_R(z) f_E(z_*) \dd z \dd z_*, \quad J = R, \\[4mm]
        & \int_{\R_+ \times \R_+} \kappa_P^{R + E \rightarrow P}(z,z_*) \mathbb{E}\big(\varphi(z')\big) f_R(z) f_E(z_*) \dd z \dd z_*, \quad J = P,
    \end{split}
\end{equation*}
and the same is true for post-interaction sizes $z'$, which now read
\begin{equation*}
    \begin{array}{ll}
        z' = z - \psi_E^{R + E \rightarrow P}(z_*) z + \eta_E^{R + E \rightarrow P}(z_*) z, & J = R, \\[4mm]
        z' = \psi_E^{R + E \rightarrow P}(z_*) z + \tilde{\eta}_E^{R + E \rightarrow P}(z_*) z, & J = P.
    \end{array}
\end{equation*}

\medskip
\noindent \textbf{Processes $\mathbf{P + P \rightarrow P}$ and $\mathbf{R + R \rightarrow R}$.} The last two processes characterize intraspecific competitive interactions that may happen among both active $P$ and hidden $R$ prey, resulting in a loss and redistribution of individuals from the same state. Therefore, only one collision operator is needed to describe the loss of prey in each state $P$ and $R$, namely
\begin{equation*}
    \begin{split}
        \int_{\R_+} \varphi(z) Q_P^{P + P \rightarrow P}(f_P,f_P)(z) \dd z &= \int_{\R_+ \times \R_+} \kappa_P^{P + P \rightarrow P}(z,z_*) \mathbb{E}\big(\varphi(z') - \varphi(z)\big) f_P(z) f_P(z_*) \dd z \dd z_*, \\[4mm]
        \int_{\R_+} \varphi(z) Q_R^{R + R \rightarrow R}(f_R,f_R)(z) \dd z &= \int_{\R_+ \times \R_+} \kappa_R^{R + R \rightarrow R}(z,z_*) \mathbb{E}\big(\varphi(z') - \varphi(z)\big) f_R(z) f_R(z_*) \dd z \dd z_*,
    \end{split}
\end{equation*}
and the post-interaction sizes $z'$ are respectively given by
\begin{equation*}
    \begin{split}
        z' &= z - \psi_P^{P + P \rightarrow P}(z_*) z + \eta_P^{P + P \rightarrow P}(z_*) z, \\[4mm]
        z' &= z - \psi_R^{R + R \rightarrow R}(z_*) z + \eta_R^{R + R \rightarrow R}(z_*) z.
    \end{split}
\end{equation*}

%%%%    SUBSECTION : MACROSCOPIC EQUATIONS    %%%%
\subsection{Moment equations for the behavioral-state densities}

In order to derive effective evolution equations for the macroscopic number densities $m_J$, $J \in \mathcal{S}$, we now multiply each equation of the Boltzmann system \eqref{eq:kinetic model} by the particular test function $\varphi(z) = z$ and integrate over $z \in \R_+$. Assuming that $z^2 f_J(t,z) \Big|_{z = 0, +\infty} = 0$ for a.e. $t \in \R_+$ and any $J \in \mathcal{S}$, the birth operators \eqref{eq:birth operators} produce the contributions
\begin{equation*}
    \int_{\R_+} z \mathcal{R}_J(f_J)(z) \dd z = -\alpha_J \Big[ z^2 f_J(z) \Big]_0^{+\infty} + \alpha_J \int_{\R_+} z f_J(z) \dd z = \alpha_J m_J, \quad J = P,R,
\end{equation*}
and, similarly, the death operators \eqref{eq:death operators} result in
\begin{equation*}
    \int_{\R_+} z \mathcal{R}_J(f_J)(z) \dd z = - \mu\, m_J, \quad J = S,H.
\end{equation*}
Next, for a generic process $I + K \rightarrow I' + K'$, with $I \neq I'$ and $K \neq K'$, the first-order moment of the Boltzmann operator $\mathcal{Q}_J^{I + K \rightarrow I' + K'}(f_I,f_K)$ is given by
\begin{equation*}
    \begin{split}
        \int_{\R_+} z \mathcal{Q}_J^{I + K \rightarrow I' + K'}&(f_I,f_K)(z) \dd z = \\[6mm]
        & \left\{ \begin{array}{lll}
        \displaystyle -\int_{\R_+ \times \R_+} \kappa_I^{I + K \rightarrow I' + K'}(z,z_*) \psi_K^{I + K \rightarrow I' + K'}(z_*) z f_I(z) f_K(z_*) \dd z \dd z_*, & \ \text{if} & J = I,\\[6mm]
        \displaystyle -\int_{\R_+ \times \R_+} \kappa_K^{I + K \rightarrow I' + K'}(z,z_*) \psi_I^{I + K \rightarrow I' + K'}(z_*) f_I(z_*) z f_K(z) \dd z \dd z_*, & \ \text{if} & J = K, \\[6mm]
        \displaystyle \int_{\R_+ \times \R_+} \kappa_{I'}^{I + K \rightarrow I' + K'}(z,z_*) \psi_K^{I + K \rightarrow I' + K'}(z_*) z f_I(z) f_K(z_*) \dd z \dd z_*, & \ \text{if} & J = I', \\[6mm]
        \displaystyle \int_{\R_+ \times \R_+} \kappa_{K'}^{I + K \rightarrow I' + K'}(z,z_*) \psi_I^{I + K \rightarrow I' + K'}(z_*) f_I(z_*) z f_K(z) \dd z \dd z_*, & \ \text{if} & J = K',
        \end{array} \right.
    \end{split}
\end{equation*}
where we have used that $\mathbb{E}\left(\eta_J^{I + K \rightarrow I' + K'}(z_*)\right) = \mathbb{E}\left(\tilde{\eta}_J^{I + K \rightarrow I' + K'}(z_*)\right) = 0$. Assuming that the transition kernels satisfy
\begin{equation*}
    \begin{split}
        \kappa_I^{I + K \rightarrow I' + K'}(z,z_*) &= \kappa_{I'}^{I + K \rightarrow I' + K'}(z,z_*), \\[4mm]
        \kappa_K^{I + K \rightarrow I' + K'}(z,z_*) &= \kappa_{K'}^{I + K \rightarrow I' + K'}(z,z_*),
    \end{split}
\end{equation*}
for a.e. $z,z_* \in \R_+$, and that there exists a common positive constant $c_{IK}^{I'K'}$ such that
\begin{equation} \label{eq:balance assumption}
    \kappa_I^{I + K \rightarrow I' + K'}(z,z_*) \psi_K^{I + K \rightarrow I' + K'}(z_*) = \kappa_K^{I + K \rightarrow I' + K'}(z,z_*) \psi_I^{I + K \rightarrow I' + K'}(z_*) = c_{IK}^{I'K'} z_*,
\end{equation}
we easily compute
\begin{equation*}
    \int_{\R_+} z \mathcal{Q}_J^{I + K \rightarrow I' + K'}(f_I,f_K)(z) \dd z = \left\{ \begin{array}{lll}
        -c_{IK}^{I'K'} m_I m_K, & \ \text{if} & J = I, K, \\[6mm]
        c_{IK}^{I'K'} m_I m_K, & \ \text{if} & J = I', K'.
        \end{array} \right.
\end{equation*}
In particular, $c_{IK}^{I'K'}$ denotes the effective transition rate of the process $I + K \rightarrow I' + K'$, determining how quickly individuals leave the original states $I,\, K$ to enter the destination ones $I',\, K'$.

\begin{remark}[Example of admissible transition kernels]
    A simple example satisfying conditions \eqref{eq:balance assumption} is obtained by considering the interaction functions $\psi_J^{I + K \rightarrow I' + K'}(z_*)$ given in Remark \ref{rem:positivity of post-interactions} and choosing interaction kernels that are proportional to the size of the interacting group \cite{BonMenTosZan}, namely
    \begin{equation*}
        \kappa_J^{I + K \rightarrow I' + K'}(z,z_*) = \kappa_J^{I + K \rightarrow I' + K'}(z_*) = \kappa_{IK}^{I'K'}(1+z_*),
    \end{equation*}
    where the positive constant $\kappa_{IK}^{I'K'}$ represents the baseline interaction frequency of process $I + K \rightarrow I' + K'$. With these choices, we obtain precisely
    \begin{equation*}
        \kappa_I^{I + K \rightarrow I' + K'}(z,z_*) \psi_K^{I + K \rightarrow I' + K'}(z_*) = \kappa_K^{I + K \rightarrow I' + K'}(z,z_*) \psi_I^{I + K \rightarrow I' + K'}(z_*) = \kappa_{IK}^{I'K'} \psi_{IK}^{I'K'} z_*,
    \end{equation*}
    and we can set $c_{IK}^{I'K'} = \kappa_{IK}^{I'K'} \psi_{IK}^{I'K'}$.
\end{remark}

If we apply these considerations to the nine processes \eqref{eq:S+P to H}--\eqref{eq:R+R to R}, we then recover
\begin{equation} \label{eq:moments of collision operators}
    \begin{split}
        & \int_{\R_+} z Q_J^{S + P \rightarrow H}(f_S,f_P)(z) \dd z = \left\{ \begin{array}{ll}
        -c_{SP}^{H} m_P m_S, & J = S, \\[4mm]
        -c_{SP}^{H} m_P m_S, & J = P, \\[4mm]
        c_{SP}^{H} m_P m_S, & J = H,
        \end{array} \right. \\[6mm]
        & \int_{\R_+} z Q_J^{S + P \rightarrow H + P}(f_S,f_P)(z) \dd z = \left\{ \begin{array}{ll}
        -c_{SP}^{HP} m_P m_S, & J = S, \\[4mm]
        c_{SP}^{HP} m_P m_S, & J = H,
        \end{array} \right. \\[6mm]
        & \int_{\R_+} z Q_S^{H + E \rightarrow H + S}(f_H,f_E)(z) \dd z = c_{HE}^{HS} m_H, \\[8mm]
        & \int_{\R_+} z Q_J^{H + E \rightarrow S}(f_H,f_E)(z) \dd z = \left\{ \begin{array}{ll}
        -c_{HE}^{S} m_H, & J = H, \\[4mm]
        c_{HE}^{S} m_H, & J = S,
        \end{array} \right. \\[6mm]
        & \int_{\R_+} z Q_J^{P + S \rightarrow R + S}(f_P,f_S)(z) \dd z = \left\{ \begin{array}{ll}
        -c_{PS}^{RS} m_P m_S, & J = P, \\[4mm]
        c_{PS}^{RS} m_P m_S, & J = R,
        \end{array} \right. \\[6mm]
        & \int_{\R_+} z Q_J^{P + H \rightarrow R + H}(f_P,f_H)(z) \dd z = \left\{ \begin{array}{ll}
        -c_{PH}^{RH} m_P m_H, & J = P, \\[4mm]
        c_{PH}^{RH} m_P m_H, & J = R,
        \end{array} \right. \\[6mm]
        & \int_{\R_+} z Q_J^{R + E \rightarrow P}(f_R,f_E)(z) \dd z = \left\{ \begin{array}{ll}
        -c_{RE}^{P} m_R, & J = R, \\[4mm]
        c_{RE}^{P} m_R, & J = P,
        \end{array} \right. \\[6mm]
        & \int_{\R_+} z Q_P^{P + P \rightarrow P}(f_P,f_P)(z) \dd z = -c_{PP}^{P} m_P^2, \\[6mm]
        & \int_{\R_+} z Q_R^{R + R \rightarrow R}(f_R,f_R)(z) \dd z = -c_{RR}^{R} m_R^2,
    \end{split}
\end{equation}
where we have used that the distribution of resources $f_E(z)$ has unit first moment $m_E = 1$.

Gathering all computations and performing simple algebraic manipulations, under the assumption \eqref{eq:balance assumption} linking the transition kernels $\kappa_J^{I + K \rightarrow I' + K'}(z,z_*)$ with the interaction functions $\psi_J^{I + K \rightarrow I' + K'}(z_*)$, the kinetic model \eqref{eq:kinetic model} reduces to an ODE system for the temporal evolution of the state densities $m_J$, $J \in \mathcal{S}$, having the form
\begin{equation} \label{eq:macroscopic BDA}
    \begin{split}
        \frac{\dd}{\dd t} m_P &= \alpha_P m_P \left( 1 - \frac{c_{PP}^{P}}{\alpha_P} m_P \right) + c_{RE}^{P} m_R - m_P \Big( \big( c_{PS}^{RS} + c_{SP}^{H} \big) m_S + c_{PH}^{RH} m_H \Big), \\[6mm]
        \frac{\dd}{\dd t} m_R &= \alpha_R m_R \left( 1 - \frac{c_{RR}^{R}}{\alpha_R} m_R \right) - c_{RE}^{P} m_R + m_P \Big( c_{PS}^{RS} m_S + c_{PH}^{RH} m_H \Big), \\[6mm]
        \frac{\dd}{\dd t} m_S &= -\mu m_S - \big( c_{SP}^{HP} + c_{SP}^{H} \big) m_P m_S + \big( c_{HE}^{S} + c_{HE}^{HS} \big) m_H, \\[6mm]
        \frac{\dd}{\dd t} m_H &= -\big( \mu + c_{HE}^{S} \big) m_H + \big( c_{SP}^{HP} + c_{SP}^{H} \big) m_P m_S.
    \end{split}
\end{equation}

\subsection{A directional velocity model for 2D spatial diffusion}

We have shown how to construct suitable kinetic operators that provide a mesoscopic formulation of predator--prey dynamics with similar reactions that lead to Beddington--De Angelis functional responses in a suitable fast-reaction asymptotics \cite{ConDesSor}. Our next goal is to extend this framework to also account for the spatial spreading of species in the two-dimensional plane, with the aim of recovering diffusion-like effects in the macroscopic system \eqref{eq:macroscopic BDA}.

The idea is to employ a 2D spatial formulation where we assume that each group of individuals of size~$z \in \R_+$ in the physical system is also characterized at the microscopic level by a spatial position $x \in \R^2$ and by a directional velocity $v \in \mathbb{S}^1$ in the plane. We then associate to each state $J \in \mathcal{S}$ a distribution function that, without risk of confusion, we still denote by $f_J = f_J(t,x,v,z)$, now depending on time $t \in \R_+$, position $x \in \R^2$, orientation $v \in \mathbb{S}^1$ and group size $z \in \R_+$. Furthermore, we introduce the corresponding relevant macroscopic observables of the state density 
\begin{equation}
    m_J(t,x) = \int_{\mathbb{S}^1 \times \R_+} z f_J(t,x,v,z) \dd v \dd z, \label{eq:macroscopic density}  
\end{equation}
and orientation flux
\begin{equation*}
    F_J(t,x) = \int_{\mathbb{S}^1 \times \R_+} v\, z f_J(t,x,v,z) \dd v \dd z, 
\end{equation*}
where we denote with $\dd v$ the normalized arc-length measure (whenever we integrate over the unit circle, we systematically divide by its total length $2\pi$). We also assume that the distribution of resources in this directional framework takes the form $f_E(v,z) = f_E(z)$ where $f_E(z)$ has unit first moment. This amounts to considering a fixed environment in velocity (a quite natural hypothesis), and will allow us to get rid once again of the density $m_E = 1$ when this will appear in the following computations of the moments.

To avoid overcomplicating the model, we maintain the simplest possible structure for the kinetic operators by assuming that the processes of predator--prey interactions and reorientation of species are modified independently by the dynamics. To start, we assume that each species is advected in straight lines along the orientations $v \in \mathbb{S}^1$. The changes in velocity are then modeled via BGK-like operators that prescribe the relaxation of each $f_J(t,x,v,z)$ toward a suitable local equilibrium state (defined below) in the tensorized form $m_J(t,x) M_J(v,z)$, which is associated with the macroscopic density $m_J$. Moreover, we assume that the birth and death processes prescribed by the redistribution operators \eqref{eq:birth operators}--\eqref{eq:death operators} and all intra-specific and inter-specific interactions described by the reactive Boltzmann operators \eqref{eq:Q_P}--\eqref{eq:Q_H} affect the change in velocity of only one state $J \in \mathcal{S}$ at a time. In other words, we require that the action of these operators in the velocity variable is at most quasi-linear with respect to the distribution of the main affected species. Given a fixed velocity $v \in \mathbb{S}^1$, the heuristics behind such modeling choices is that the birth of prey (active or hidden) and the death of predators (searching or handling) should occur along the direction of movement $v$ under consideration, and any passage of individuals from one state to another following processes \eqref{eq:S+P to H}--\eqref{eq:R+R to R} should not alter the orientation $v$, whatever the velocity $v_* \in \mathbb{S}^1$ of the other state involved in the microscopic interaction.

Based on these considerations, we model the evolution of the vector distribution $\mathbf{f} = (f_P,f_R,f_S,f_H,f_E)$ through the kinetic system
\begin{equation} \label{eq:kinetic model new}
    \frac{\partial}{\partial t} f_J + \bar{v}_J v \cdot \nabla_x f_J = \mathcal{L}_J(f_J) + \mathcal{R}_J(f_J) + \mathcal{Q}_J(\mathbf{f},\mathbf{f}), \quad t \in \R_+, \, x \in \R^2, \, v \in \mathbb{S}^1, \, z \in \R_+, \, J \in \mathcal{S},
\end{equation}
where the left-hand side describes the free transport of species, with $\bar{v}_J > 0$ denoting the characteristic velocity of individuals belonging to state $J \in \mathcal{S}$, while the right-hand side models all kinetic-type interactions in $v \in \mathbb{S}^1$ and $z \in \R_+$. More precisely, changes in orientation are modeled as linear relaxation operators~$\mathcal{L}_J(f_J),\,J \in \mathcal{S}$, of the form \cite{BGK}
\begin{equation*} 
    \mathcal{L}_J(f_J)(t,x,v,z) = \nu_J \big( m_J(t,x) M_J(v,z) - f_J(t,x,v,z) \big),
\end{equation*}
which prescribe the relaxation of each $f_J$ toward a specific distribution $m_J M_J$, at a constant characteristic relaxation speed $\nu_J > 0$. Here, we say that the vector distribution function $\mathbf{M} = (M_J)_{J \in \mathcal{S}}$ defines an equilibrium state (in velocity and size) for the system if each component $M_J = M_J(v,z)$ is nonnegative~$M_J(v,z) \geq 0$ for a.e. $v \in \mathbb{S}^1$ and $z \in \R_+$, satisfies the minimal regularity $M_J \in L^1(\mathbb{S}^1 \times \R_+)$, normalizes the number of individuals in state $J$, namely
\begin{equation} \label{eq:property 1 of equilibria}
    \int_{\mathbb{S}^1 \times \R_+} z M_J(v,z) \dd v \dd z = 1,
\end{equation}
and has zero flux, i.e.
\begin{equation} \label{eq:property 2 of equilibria}
    \int_{\mathbb{S}^1 \times \R_+} v\, z M_J(v,z) \dd v \dd z = 0.
\end{equation}
Notice in particular that by construction the operators $\mathcal{L}_J(f_J)$, $J \in \mathcal{S}$, satisfy the natural conservation laws
\begin{equation*} 
    \int_{\mathbb{S}^1 \times \R_+} z \big(m_J M_J(v,z) - f_J(v,z) \big) \dd v \dd z = 0,
\end{equation*}
thanks to the properties \eqref{eq:property 1 of equilibria} and \eqref{eq:property 2 of equilibria} of the local equilibria $M_J$. We point out that one could also resort to more elaborate reorientation models, like the classical kinetic tumbling operators which find application in other biological contexts, e.g. the motion of cells and bacteria \cite{Alt,OthDunAlt,Str}. However, the use of BGK-type operators simplifies the exposition and imposes a tensorized structure on the relaxation states  $m_J M_J$, which ultimately allows us to provide explicit examples of ecologically meaningful equilibrium distributions (see Remark \ref{rem:equilibria} below).

On the other hand, for the kinetic updates in the size $z \in \R_+$ we still consider the same birth and death operators \eqref{eq:birth operators} and \eqref{eq:death operators}, while we assume that the action of the Boltzmann operators $\mathcal{Q}_J(\mathbf{f},\mathbf{f})$, $J \in \mathcal{S}$, is now modulated by a further integration over the velocity domain, which corresponds to characterizing any component $\mathcal{Q}_J^{I + K \rightarrow I' + K'}(f_I,f_K)$ in weak form as
\begin{equation*}
    \begin{split}
        \int_{\R_+} \varphi(z) & \mathcal{Q}_J^{I + K \rightarrow I' + K'}(f_I,f_K)(v,z) \dd z = \\[2mm]
        & \int_{\mathbb{S}^1 \times \R_+ \times \R_+} \kappa_J^{I + K \rightarrow I' + K'}(z,z_*) \mathbb{E}\big(\varphi(z')-\delta_{JI}\varphi(z)\big) f_I(v,z) f_K(v_*,z_*) \dd v_* \dd z \dd z_*, \quad J = I,I', \\[4mm]
        & \int_{\mathbb{S}^1 \times \R_+ \times \R_+} \kappa_J^{I + K \rightarrow I' + K'}(z,z_*) \mathbb{E}\big(\varphi(z')-\delta_{JK}\varphi(z)\big) f_I(v_*,z_*) f_K(v,z) \dd v_* \dd z \dd z_*, \quad J = K,K',
    \end{split}
\end{equation*}
where the post-interaction size $z'$ is still determined through relations \eqref{eq:microsopic interactions}. This simply means that we average the interactions $(z,z_*) \mapsto z'$ over all possible orientations $v_* \in \mathbb{S}^1$.

\begin{remark}[Examples of equilibrium states] \label{rem:equilibria}
    The class of equilibria $\mathbf{M}$ is obviously non-empty and contains both isotropic and non-isotropic distributions in $v$. The most interesting examples in the isotropic case (distributions that are constant along the orientation variable) are represented by uniform distributions on bounded domains $z \in [0,Z_J]$, namely
    \begin{equation*}
        M_J(v,z) = \frac{2}{Z_J^2} \chi_{[0,Z_J]},
    \end{equation*}
    and generalized Gamma densities of the form
    \begin{equation*}
        M_J(v,z) = \frac{r_J}{b_J^{a_J+1} \Gamma\left( \frac{a_J+1}{r_J} \right)} z^{a_J-1} e^{-(z/b_J)^{r_J}},
    \end{equation*}
    for some shape and scale parameters $a_J, \, b_J, \, r_J > 0$, where $\Gamma(s)$, $s > 0$, denotes the standard Gamma function. The latter equilibria are particularly relevant, since the long-time emergence of Gamma distributions was recently investigated in the context of solutions to Fokker--Planck-type models for Lotka--Volterra interactions~\cite{BonMenTosZan,TosZan_eq}.
\end{remark}

%%%%%%%%%%%%%%%%%%%%%%%%%%%%%%%%%%%%%%%%%%
%%%%    SECTION 3: DIFFUSIVE LIMIT    %%%%
%%%%%%%%%%%%%%%%%%%%%%%%%%%%%%%%%%%%%%%%%%

\section{Diffusive fast-reaction asymptotics of the kinetic model}\label{sec:diff}

Starting from the kinetic system \eqref{eq:kinetic model new}, we now formally derive the macroscopic reaction--cross-diffusion system with a Beddington--DeAngelis functional response obtained  in~\cite{DesSor} by combining a diffusive scaling with a fast-reaction regime. For this, we assume that the relaxation of each $f_J$ toward the prescribed local equilibrium state~$M_J$ is the dominant dynamics. Then, fast reactions constitute the next lower-order processes, occurring at the same scale as the spatial transport of species. At last, the lowest-order processes are the slow reactions, occurring at the same scale as the temporal derivative. More precisely, let us introduce a small parameter~$\eps \ll 1$ and a parabolic rescaling $(t,x) \mapsto (t/\eps^2,x/\eps)$ of time and space. Accordingly, we define the vector distribution $\mathbf{f}^\eps = (f_P^\eps,f_R^\eps,f_S^\eps,f_H^\eps,f_E)$ in the rescaled variables as $f_J^\eps = f_J(t/\eps^2,x/\eps,v,z)$, together with the corresponding macroscopic densities and fluxes
\begin{equation*}
	\begin{split}
		m_J^\eps(t,x) &= \int_{\mathbb{S}^1 \times \R_+} z f_J^\eps(t,x,v,z) \dd v \dd z, \\[4mm]
		F_J^\eps(t,x) &= \int_{\mathbb{S}^1 \times \R_+} v \, z f_J^\eps(t,x,v,z) \dd v \dd z,
	\end{split}
\end{equation*}
for any $J \in \mathcal{S}$. We highlight in particular that the distribution of resources $f_E$ does not depend on $\eps$. Then, the rescaled version of equations \eqref{eq:kinetic model new}, modeling the evolution of $\mathbf{f}^\eps$, is given by the diffusive system
\begin{equation} \label{eq:rescaled kinetic model}
    \eps \frac{\partial}{\partial t} f_J^\eps + \bar{v}_J v \cdot \nabla_x f_J^\eps = \frac{1}{\eps} \big( \mathcal{L}_J(f_J^\eps) + \mathcal{R}_J^\eps(f_J^\eps) + \mathcal{Q}_J^\eps(\mathbf{f}^\eps,\mathbf{f}^\eps) \big), \quad t \in \R_+, \, x \in \R^2, \, v \in \mathbb{S}^1, \, z \in \R_+, \, J \in \mathcal{S},
\end{equation}
where now the BGK operators
\begin{equation*}
	\mathcal{L}_J(f_J^\eps)(t,x,v,z) = \nu_J(m_J^\eps(t,x) M_J(v,z) - f_J^\eps(t,x,v,z)), \quad J \in \mathcal{S},
\end{equation*}
prescribe the relaxation toward $\eps$-dependent tensorized local equilibria  $m_J^\eps M_J$ (notice that the state~$\mathbf{M}$ remains untouched by the rescaling), while the rescaled redistribution operators $\mathcal{R}_J^\eps(f_J^\eps)$ and the Boltzmann operators $\mathcal{Q}_J^\eps(\mathbf{f}^\eps,\mathbf{f}^\eps)$ depend themselves on $\eps$ through the slow--fast formulations
\begin{equation*}
    \begin{split}
        \mathcal{R}_J^\eps(f_J^\eps) &= \eps^2 \mathcal{R}_J(f_J^\eps), \quad J \in \mathcal{S}, \\[6mm]
        \mathcal{Q}_P^\eps(\mathbf{f}^\eps,\mathbf{f}^\eps) &= \eps^2 \mathcal{Q}_P^{S + P \rightarrow H}(f_S^\eps,f_P^\eps) + \eps \mathcal{Q}_P^{P + S \rightarrow R + S}(f_P^\eps,f_S^\eps) + \eps \mathcal{Q}_P^{P + H \rightarrow R + H}(f_P^\eps,f_H^\eps) \\[4mm]
        &\qquad + \eps \mathcal{Q}_P^{R + E \rightarrow P}(f_R^\eps,f_E) + \eps^2 \mathcal{Q}_P^{P + P \rightarrow P}(f_P^\eps,f_P^\eps), \\[6mm]
        \mathcal{Q}_R^\eps(\mathbf{f}^\eps,\mathbf{f}^\eps) &= \eps \mathcal{Q}_R^{P + S \rightarrow R + S}(f_P^\eps,f_S^\eps) + \eps \mathcal{Q}_R^{P + H \rightarrow R + H}(f_P^\eps,f_H^\eps) + \eps \mathcal{Q}_R^{R + E \rightarrow P}(f_R^\eps,f_E) + \eps^2 \mathcal{Q}_R^{R + R \rightarrow R}(f_R^\eps,f_R^\eps), \\[6mm]
        \mathcal{Q}_S^\eps(\mathbf{f}^\eps,\mathbf{f}^\eps) &= \eps^2 \mathcal{Q}_S^{S + P \rightarrow H}(f_S^\eps,f_P^\eps) + \eps \mathcal{Q}_S^{S + P \rightarrow H + P}(f_S^\eps,f_P^\eps) + \eps^2 \mathcal{Q}_S^{H + E \rightarrow H + S}(f_H^\eps,f_E) + \eps \mathcal{Q}_S^{H + E \rightarrow S}(f_H^\eps,f_E), \\[6mm]
        \mathcal{Q}_H^\eps(\mathbf{f}^\eps,\mathbf{f}^\eps) &= \eps^2 \mathcal{Q}_H^{S + P \rightarrow H}(f_S^\eps,f_P^\eps) + \eps \mathcal{Q}_H^{S + P \rightarrow H + P}(f_S^\eps,f_P^\eps) + \eps \mathcal{Q}_H^{H + E \rightarrow S}(f_H^\eps,f_E),
    \end{split}
\end{equation*}
based on the following slow--fast microscopic dynamics for predators
\begin{equation*}
    \begin{array}{lll}
        S + P \rightarrow H, & & \textrm{(slow process)} \\[2mm]
        S + P \rightarrow H + P, & & \textrm{(fast process)}  \\[2mm]
        H + E \rightarrow H + S, & & \textrm{(slow process)}  \\[2mm]
        H + E \rightarrow S. & & \textrm{(fast process)} 
    \end{array}
\end{equation*}
and for prey
\begin{equation*}
    \begin{array}{lll}
        P + S \rightarrow R + S, & & \textrm{(fast process)} \\[2mm]
        P + H \rightarrow R + H, & & \textrm{(fast process)} \\[2mm]
        R + E \rightarrow P, & & \textrm{(fast process)} \\[2mm]
        P + P \rightarrow P, & & \textrm{(slow process)} \\[2mm]
        R + R \rightarrow R. & & \textrm{(slow process)}
    \end{array}
\end{equation*}
In order to investigate the formal diffusion asymptotics $\eps \ll 1$ of model \eqref{eq:rescaled kinetic model}, we employ the moment method to characterize evolution equations (at an intermediate level between the kinetic and macroscopic descriptions) for the number densities $m_J^\eps$ and for the fluxes $F_J^\eps$, $J \in \mathcal{S}$, to recover in the limit $\eps \to 0$ a macroscopic reaction--cross-diffusion model like the one proposed in \cite{DesSor}. Initially, we notice that the dominant order terms of system \eqref{eq:rescaled kinetic model} are the BGK operators which drive the relaxation processes in $v$. This suggests that the relations $\mathcal{L}_J(f_J^\eps) = \mathcal{O}(\eps)$, $J \in \mathcal{S}$, hold at leading order, allowing to infer the approximations~$f_J^\eps \underset{\eps \to 0}{\sim} m_J^\eps M_J$ for any $J \in \mathcal{S}$. Therefore, it is natural to expand $\mathbf{f}^\eps$ as
\begin{equation} \label{eq:ansatz}
    f_J^\eps(t,x,v,z) = m_J^\eps(t,x) M_J(v,z) + \eps \tilde{f}_J^\eps(t,x,v,z), \quad J \in \mathcal{S},
\end{equation}
where the perturbations $\tilde{f}_J^\eps$ satisfy, by construction (since the density $m_J^\eps$ of each state is solely carried by the equilibrium $m_J^\eps M_J$), the zero-mass conditions
\begin{equation} \label{eq:property of perturbations}
    \int_{\mathbb{S}^1 \times \R_+} z \tilde{f}_J^\eps(v,z) \dd v \dd z = 0,
\end{equation}
and are assumed to be bounded in some suitable normed space. In particular, since the distribution of resources $f_E$ is already at equilibrium, we also set $M_E = f_E$ and $\tilde{f}_E = 0$. We can thus plug the Ansatz~\eqref{eq:ansatz} into system~\eqref{eq:rescaled kinetic model} to recover evolution equations for the vector distribution $\mathbf{\tilde{f}}^\eps = (\tilde{f}_P^\eps,\tilde{f}_R^\eps,\tilde{f}_S^\eps,\tilde{f}_H^\eps,0)$, depending on the local equilibrium state $\mathbf{m^\eps M} = (m_P^\eps M_P,m_R^\eps M_R,m_S^\eps M_S,m_H^\eps M_H,f_E)$. These are given, for any $J \in \mathcal{S}$,~by
\begin{equation} \label{eq:evolution of perturbations}
    \begin{split}
        \left(\eps \frac{\partial}{\partial t} m_J^\eps + \bar{v}_J v \cdot \nabla_x m_J^\eps \right) M_J & + \eps^2 \frac{\partial}{\partial t} \tilde{f}_J^\eps + \eps \bar{v}_J v \cdot \nabla_x \tilde{f}_J^\eps = -\nu_J \tilde{f}_J^\eps + \frac{1}{\eps} \mathcal{R}_J^\eps(m_J^\eps M_J) + \mathcal{R}_J^\eps(\tilde{f}_J^\eps) \\[4mm]
        & + \frac{1}{\eps} \mathcal{Q}_J^\eps(\mathbf{m^\eps M},\mathbf{m^\eps M}) +\mathcal{Q}_J^\eps(\mathbf{m^\eps M},\mathbf{\tilde{f}}^\eps) + \mathcal{Q}_J^\eps(\mathbf{\tilde{f}}^\eps,\mathbf{m^\eps M}) + \eps \mathcal{Q}_J^\eps(\mathbf{\tilde{f}}^\eps,\mathbf{\tilde{f}}^\eps),
    \end{split}
\end{equation}
where we have used that $\mathcal{L}_J(m_J^\eps M_J) = 0$.

Our aim is now to apply the moment method to the above equation, as a means to identify its macroscopic leading-order terms in $\eps$. To do this, we compute the evolution of the zeroth-order moments in velocity by multiplying each equation in the system \eqref{eq:evolution of perturbations} by $z$ and integrating over $\mathbb{S}^1 \times \R_+$. Using that each $M_J$ has unitary macroscopic density \eqref{eq:property 1 of equilibria} and zero flux \eqref{eq:property 2 of equilibria}, together with the zero-mass condition \eqref{eq:property of perturbations} satisfied by the perturbations $\tilde{f}_J^\eps$, it is easy to check that
\begin{equation*}
    \begin{split}
        & \int_{\mathbb{S}^1 \times \R_+} z \left(\eps \frac{\partial}{\partial t} m_J^\eps + \bar{v}_J v \cdot \nabla_x m_J^\eps \right) M_J(v,z) \dd v \dd z = \eps \frac{\partial}{\partial t} m_J^\eps, \hspace{6cm}
    \end{split}
\end{equation*}
\begin{equation*}
    \begin{split}
        & \int_{\mathbb{S}^1 \times \R_+} z \mathcal{R}_J^\eps(m_J^\eps M_J)(v,z) \dd v \dd z = -\eps^2 \alpha_J m_J^\eps \int_{\mathbb{S}^1} \Big[ z^2 M_J(v,z) \Big]_0^{+\infty} \dd v + \eps^2 \alpha_J m_J^\eps \int_{\mathbb{S}^1 \times \R_+} z M_J(z) \dd z \\[3mm]
        & \hspace{5.15cm} = \eps^2 \alpha_J m_J^\eps, \quad J = P,R, \\[4mm]
        & \int_{\mathbb{S}^1 \times \R_+} z \mathcal{R}_J^\eps(m_J^\eps M_J)(v,z) \dd v \dd z = -\eps^2 \mu m_J^\eps, \quad J = S,H,
    \end{split}
\end{equation*}
assuming that the boundary conditions $z^2 M_J(t,x,v,z) \Big|_{z = 0,+\infty} = 0$ hold for a.e. $t \in \R_+$, $x \in \R^2$, $v \in \mathbb{S}^1$ and any $J \in \mathcal{S}$. Moreover, given a generic transition $I + K \rightarrow I' + K'$, simple computations show that any integral term involving $Q_J^{I + K \rightarrow I' + K'}(m_I^\eps M_I,m_K^\eps M_K)$ produces a contribution
\begin{equation*}
    \begin{split}
        \int_{\mathbb{S}^1 \times \R_+} z \mathcal{Q}_J^{I + K \rightarrow I' + K'}(m_I^\eps M_I,m_K^\eps M_K)(v,z) \dd v \dd z
        = \left\{ \begin{array}{lll}
        -c_{IK}^{I'K'} m_I^\eps m_K^\eps, & \ \text{if} & J = I, K, \\[6mm]
        c_{IK}^{I'K'} m_I^\eps m_K^\eps, & \ \text{if} & J = I', K',
        \end{array} \right.
    \end{split}
\end{equation*}
while the first-order moment of all other reactive terms $Q_J^{I + K \rightarrow I' + K'}(m_I^\eps M_I,\tilde{f}_K^\eps)$ and $Q_J^{I + K \rightarrow I' + K'}(\tilde{f}_I^\eps,\tilde{f}_K^\eps)$ vanishes
\begin{equation*}
    \begin{split}
        & \int_{\mathbb{S}^1 \times \R_+} z \mathcal{Q}_J^{I + K \rightarrow I' + K'}(m_I^\eps M_I,\tilde{f}_K^\eps)(v,z) \dd v \dd z = 0, \\[4mm]
        & \int_{\mathbb{S}^1 \times \R_+} z \mathcal{Q}_J^{I + K \rightarrow I' + K'}(\tilde{f}_I^\eps,\tilde{f}_K^\eps)(v,z) \dd v \dd z = 0,
    \end{split}
\end{equation*}
since the perturbations $\tilde{f}_J^\eps$ have zero mass \eqref{eq:property of perturbations}. Therefore, taking into account these considerations and repeating the computations made to derive the moments \eqref{eq:moments of collision operators}, we find that equations \eqref{eq:evolution of perturbations} reduce to a fast-reaction-like system prescribing the evolution of the macroscopic densities $\mathbf{m}^\eps = (m_J^\eps)_{J \in \mathcal{S}}$, given by
\begin{equation} \label{eq:zeroth-order moments}
    \begin{split}
        &\frac{\partial}{\partial t} m_P^\eps + \bar{v}_P \nabla_x \cdot \tilde{F}_P^\eps = \alpha_P m_P^\eps \left( 1 - \frac{c_{PP}^{P}}{\alpha_P} m_P^\eps \right) - \frac{1}{\eps} m_P^\eps \Big( c_{PS}^{RS} m_S^\eps + c_{PH}^{RH} m_H^\eps \Big) + \frac{c_{RE}^{P}}{\eps} m_R^\eps - c_{SP}^{H} m_P^\eps m_S^\eps, \\[6mm]
        &\frac{\partial}{\partial t} m_R^\eps + \bar{v}_R \nabla_x \cdot \tilde{F}_R^\eps = \alpha_R m_R^\eps \left( 1 - \frac{c_{RR}^{R}}{\alpha_R} m_R^\eps \right) + \frac{1}{\eps} m_P^\eps \Big( c_{PS}^{RS} m_S^\eps + c_{PH}^{RH} m_H^\eps \Big) - \frac{c_{RE}^{P}}{\eps} m_R^\eps, \\[6mm]
        &\frac{\partial}{\partial t} m_S^\eps + \bar{v}_S \nabla_x \cdot \tilde{F}_S^\eps = -\mu m_S^\eps + \frac{1}{\eps} c_{HE}^{S} m_H^\eps - \frac{1}{\eps} c_{SP}^{HP} m_P^\eps m_S^\eps - c_{SP}^{H} m_P^\eps m_S^\eps + c_{HE}^{HS} m_H^\eps, \\[6mm]
        & \frac{\partial}{\partial t} m_H^\eps + \bar{v}_H \nabla_x \cdot \tilde{F}_H^\eps = -\mu m_H^\eps - \frac{1}{\eps} c_{HE}^{S} m_H^\eps + \frac{1}{\eps} c_{SP}^{HP} m_P^\eps m_S^\eps + c_{SP}^{H} m_P^\eps m_S^\eps,
    \end{split}
\end{equation}
where we have denoted the flux of each perturbation $\tilde{f}_J^\eps$ as
\begin{equation*}
    \tilde{F}_J^\eps(t,x) = \int_{\mathbb{S}^1 \times \R_+} v\, z \tilde{f}_J^\eps(t,x,v,z) \dd v \dd z.
\end{equation*}
In order to close system \eqref{eq:zeroth-order moments}, we need to determine a relation between the densities $m_J^\eps$ and the fluxes $\tilde{F}_J^\eps$ that should link them together in the limit $\eps \to 0$. We then proceed by multiplying system \eqref{eq:evolution of perturbations} by $v\, z$ and integrating over $\mathbb{S}^1 \times \R_+$, to recover the equations governing the evolution of the first-order moments, namely
\begin{equation} \label{eq:first-order moments}
    \begin{split}
        \eps^2 \frac{\partial}{\partial t} \tilde{F}_J^\eps & + \bar{v}_J \nabla_x m_J^\eps \cdot \int_{\mathbb{S}^1 \times \R_+} v \otimes v \,z  M_J(v,z) \dd v \dd z + \eps \bar{v}_J \nabla_x \cdot \int_{\mathbb{S}^1 \times \R_+} v \otimes v \,z  \tilde{f}_J^\eps(v,z) \dd v \dd z = -\nu_J \tilde{F}_J^\eps \\[4mm]
        & + \frac{1}{\eps} \int_{\mathbb{S}^1 \times \R_+} v\, z \mathcal{R}_J^\eps(m_J^\eps M_J) \dd v \dd z + \int_{\mathbb{S}^1 \times \R_+} v\, z \mathcal{R}_J^\eps(\tilde{f}_J^\eps) \dd v \dd z + \frac{1}{\eps} \int_{\mathbb{S}^1 \times \R_+} v\, z \mathcal{Q}_J^\eps(\mathbf{m^\eps M},\mathbf{m^\eps M}) \dd v \dd z \\[4mm]
        & + \int_{\mathbb{S}^1 \times \R_+} v\, z \big( \mathcal{Q}_J^\eps(\mathbf{m^\eps M},\mathbf{\tilde{f}}^\eps) + \mathcal{Q}_J^\eps(\mathbf{\tilde{f}}^\eps,\mathbf{m^\eps M}) \big) \dd v \dd z + \eps \int_{\mathbb{S}^1 \times \R_+} v\, z \mathcal{Q}_J^\eps(\mathbf{\tilde{f}}^\eps,\mathbf{\tilde{f}}^\eps) \dd v \dd z,
    \end{split}
\end{equation}
where we have used that the equilibria $M_J$ satisfy the zero-flux condition \eqref{eq:property 2 of equilibria}. This latter assumption also allows us to show that, for a generic process $I + K \rightarrow I' + K'$, the leading-order integral terms involving the operators $Q_J^{I + K \rightarrow I' + K'}(m_I^\eps M_I,m_K^\eps M_K)$ vanish, since
\begin{equation*}
    \begin{split}
        \int_{\mathbb{S}^1 \times \R_+} v \,z \mathcal{Q}_J^{I + K \rightarrow I' + K'}&(m_I^\eps M_I,m_K^\eps M_K)(v,z) \dd v \dd z = \\[6mm]
        & \left\{ \begin{array}{lll}
        \displaystyle -c_{IK}^{I'K'} m_I^\eps m_K^\eps \int_{(\mathbb{S}^1 \times \R_+)^2} v \, z M_I(v,z) z_* M_K(v_*,z_*) \dd v \dd v_* \dd z \dd z_*, & \ \text{if} & J = I, K, \\[6mm]
        \displaystyle c_{IK}^{I'K'} m_I^\eps m_K^\eps \int_{(\mathbb{S}^1 \times \R_+)^2} v \, z_* M_I(v_*,z_*) z M_K(v,z) \dd v \dd v_* \dd z \dd z_*, & \ \text{if} & J = I', K',
        \end{array} \right. \\[8mm]
        & = 0.
    \end{split}
\end{equation*}
In particular, these cancellations allow us to explicitly characterize the dominant terms (advection of~$m_J^\eps$ and flux $\tilde{F}_J^\eps$) in system \eqref{eq:first-order moments} and infer the first-order approximations
\begin{equation} \label{eq:Fick law}
    \tilde{F}_J^\eps = - C_J \nabla_x m_J^\eps + \mathcal{O}(\eps), \quad J \in \mathcal{S},
\end{equation}
where the coefficient matrices of these Fick-type laws are given by
\begin{equation} \label{eq:diffusion coefficients}
    C_J = \frac{\bar{v}_J}{\nu_J} \int_{\mathbb{S}^1 \times \R_+} v \otimes v \, z M_J(v,z) \dd v \dd z,
\end{equation}
depending on the characteristic speed $\bar{v}_J$ of the species, the rate $\nu_J$ of relaxation to equilibrium, and the second-order tensor $\int_{\mathbb{S}^1 \times \R_+} v \otimes v \, z M_J(v,z) \dd v \dd z$,  which describes how the macroscopic density $m_J^\eps$ is distributed at equilibrium across directions.

The aim is now to perform the hydrodynamic limit $\eps \to 0$ of equations \eqref{eq:zeroth-order moments}--\eqref{eq:Fick law}, from which we expect to recover a reaction--cross-diffusion system governing the evolution of the macroscopic quantities $m_P^\eps+m_R^\eps$ and~$m_S^\eps+m_H^\eps$. We begin by making the modeling assumption that the switching rate of prey from state $P$ to~$R$ when they encounter searching (process $P + S \rightarrow R + S$) or handling predators (process~$P + H \rightarrow R + H$) is the same.  Mathematically, this hypothesis amounts to fix the constant $c_{PH}^{RH} = c_{PS}^{RS}$ in system \eqref{eq:zeroth-order moments}, allowing to simplify the expression of the fast-reaction terms
\begin{equation*}
    c_{PS}^{RS} m_S^\eps + c_{PH}^{RH} m_H^\eps = c_{PS}^{RS} (m_S^\eps + m_H^\eps).
\end{equation*}
This reduction will be essential in the following derivation, in order to perform the fast-reaction limit in the right-hand side of \eqref{eq:zeroth-order moments}. The assumption means that prey respond to the mere presence of a predator regardless of whether it is searching for prey or handling another prey individual, and that shelter entry is an immediate antipredator response that does not depend strongly on the predator's current activity. It is thus reasonable from an ecological standpoint, since our model is not intended to capture more complex behaviors associated with prey perception and decision-making processes.

Our derivation is based on formal arguments, assuming that for each $J \in \mathcal{S}$ there exist a limit macroscopic density $m_J = m_J(t,x)$ and a limit flux $\tilde{F}_J = \tilde{F}_J(t,x)$ such that respectively $\displaystyle \lim_{\eps \to 0} m_J^\eps = m_J$ and $\displaystyle \lim_{\eps \to 0} \tilde{F}_J^\eps = \tilde{F}_J$. At first, one finds that when $\eps \to 0$ equations \eqref{eq:Fick law} lead to the asymptotic relations
\begin{equation} \label{eq:limit-flux}
    \tilde{F}_J = -C_J \nabla_x m_J,
\end{equation}
linking together limit fluxes and densities. Therefore, in order to determine the closed evolution of~$\mathbf{m} = (m_J)_{J \in \mathcal{S}}$, we are solely left with resolving the fast-reaction limit that appears on the right-hand side of system~\eqref{eq:zeroth-order moments}, the procedure being identical to the standard one used to deal with fast-reaction--diffusion systems~\cite{GerGyl,ConDesSor,DesSor}. For this, let us introduce the total density of prey~$m_X = m_P + m_R$ and predators~$m_Y = m_S + m_H$. To start, the fast-reaction terms imply that
\begin{align*}
    & m_R - \frac{c_{PS}^{RS}}{c_{RE}^{P}} m_P m_Y = \mathcal{O}(\eps), \\[4mm]
    & m_H - \frac{c_{SP}^{HP}}{c_{HE}^{S}} m_P m_S = \mathcal{O}(\eps),
\end{align*}
which ensure in the limit $\eps \to 0$ the recovery of the two fundamental closure relations
\begin{align}
    & m_R = \beta m_P m_Y, \label{eq:fast-reaction 1} \\[4mm]
    & m_H = \delta m_P m_S, \label{eq:fast-reaction 2}
\end{align}
where we have denoted with $\beta = c_{PS}^{RS}/c_{RE}^{P}$ and $\delta = c_{SP}^{HP}/c_{HE}^{S}$. At fixed densities, increasing $\beta$ favors the hidden state, whereas increasing $\delta$ favors the handling state. More precisely, since the closure relations give $m_R/m_P = \beta m_Y$ and $m_H/m_S = \delta m_P$, hidden prey outnumber active prey when $\beta m_Y > 1$,while handling predators outnumber searching predators when $\delta m_P > 1$.

The aim is now to express all state densities $m_J$, $J \in \mathcal{S}$, in terms of the total quantities $m_X$ and $m_Y$.  By adding and subtracting $m_P$ from relation \eqref{eq:fast-reaction 1}, we first infer that
\begin{equation} \label{eq:m_P}
    m_P = \frac{m_X}{1 + \beta m_Y},
\end{equation}
and proceeding similarly with relation \eqref{eq:fast-reaction 2} by adding and subtracting $m_S$ this time, we also get
\begin{equation} \label{eq:m_S}
    m_S = \frac{m_Y}{1 + \delta m_P} = \frac{m_Y (1 + \beta m_Y)}{1 + \delta m_X + \beta m_Y},
\end{equation}
where the second equality comes from substituting $m_P$ from \eqref{eq:m_P}. But now we may also insert \eqref{eq:m_P} into relation \eqref{eq:fast-reaction 1} to find
\begin{equation} \label{eq:m_R}
    m_R = \frac{\beta m_X m_Y}{1 + \beta m_Y},
\end{equation}
and finally expressions \eqref{eq:m_P}--\eqref{eq:m_S} allow us to replace $m_P$ and $m_S$ in relation \eqref{eq:fast-reaction 2}, obtaining
\begin{equation} \label{eq:m_H}
    m_H = \frac{\delta m_X m_Y}{1 + \delta m_X + \beta m_Y}.
\end{equation}
We now add the first and second equations of system \eqref{eq:zeroth-order moments} to eliminate the fast-reaction terms, and then we take the limit $\eps \to 0$ to recover
\begin{equation*}
    \frac{\partial}{\partial t} m_X + \nabla_x \cdot \big( \bar{v}_P \tilde{F}_P + \bar{v}_R \tilde{F}_R \big) = \alpha_P m_P \left( 1 - \frac{m_P}{K_P} \right) +  \alpha_R m_R \left( 1 - \frac{m_R}{K_R} \right) - c_{SP}^{H} m_P m_S,
\end{equation*}
where we have denoted with $K_P = \alpha_P/c_{PP}^{P}$ and $K_R = \alpha_R/c_{RR}^{R}$ the carrying capacities of states $P$ and $R$. To conclude, using Fick's laws \eqref{eq:limit-flux} to express the fluxes $\tilde{F}_P$ and~$\tilde{F}_R$ in terms of the corresponding state densities~$m_P$ and~$m_R$, we can combine them with relations~\eqref{eq:m_P}--\eqref{eq:m_R} to derive a reaction--cross-diffusion equation for~$m_X$, of the form
\begin{equation} \label{eq:limit 1}
    \begin{split}
        \frac{\partial}{\partial t} m_X - \nabla_x \cdot \Bigg( D_P \nabla_x \frac{m_X}{1 + \beta m_Y} + D_R & \nabla_x \frac{\beta m_X m_Y}{1 + \beta m_Y} \Bigg) = \alpha_P \frac{m_X}{1 + \beta m_Y} \left( 1 - \frac{c_{PP}^{P}}{\alpha_P} \frac{m_X}{1 + \beta m_Y} \right) \\[6mm]
        & + \alpha_R \frac{\beta m_X m_Y}{1 + \beta m_Y} \left( 1 - \frac{c_{RR}^{R}}{\alpha_R} \frac{\beta m_X m_Y}{1 + \beta m_Y} \right) - c_{SP}^{H} \frac{m_X m_Y}{1 + \delta m_X + \beta m_Y},
    \end{split}
\end{equation}
where we have defined the diffusion matrices $D_J = \bar{v}_J C_J$, with $C_J$ given by \eqref{eq:diffusion coefficients}.

\begin{remark}[A predator-dependent carrying capacity]
    We highlight that under the particular assumption that the intrinsic birth rates $\alpha_J$ and carrying capacities $K_J$, $J = P,\,R$, simplify to
    \begin{equation} \label{eq:simplified setting}
        \begin{split}
            \alpha_P &= \alpha_R =: \alpha, \\[4mm]
            K_P &= K_R =: K,
        \end{split}
    \end{equation}
    the corresponding growth terms in equation \eqref{eq:limit 1} reduce to
    \begin{equation*}
        \alpha_P \frac{m_X}{1 + \beta m_Y} \left( 1 - \frac{c_{PP}^{P}}{\alpha_P} \frac{m_X}{1 + \beta m_Y} \right) + \alpha_R \frac{\beta m_X m_Y}{1 + \beta m_Y} \left( 1 - \frac{c_{RR}^{R}}{\alpha_R} \frac{\beta m_X m_Y}{1 + \beta m_Y} \right) = \alpha m_X \left( 1 - \frac{m_X}{f_K(m_Y)} \right),
    \end{equation*}
    with a novel predator-dependent carrying capacity $f_K(m_Y)$ defined by
    \begin{equation*}
        f_K(m_Y) = K \frac{\big(1 + \beta m_Y\big)^2}{1 + \big(\beta m_Y\big)^2}.
    \end{equation*}
    The effect of this function and how it compares to existing models in the literature will be discussed in the next section.
\end{remark}

Starting from the addition of the third and fourth equations of system \eqref{eq:zeroth-order moments}, we can then perform similar computations to formally derive in the limit $\eps \to 0$ an evolution equation for $m_Y$, which reads
\begin{equation*}
    \frac{\partial}{\partial t} m_Y + \nabla_x \cdot \big( \bar{v}_S \tilde{F}_S + \bar{v}_H \tilde{F}_H \big) = -\mu m_Y + c_{HE}^{HS} m_H.
\end{equation*}
Using once again the asymptotic relations \eqref{eq:limit-flux} on the fluxes, together with the explicit expressions \eqref{eq:m_S}--\eqref{eq:m_H} of $m_S$ and $m_H$, we finally recover
\begin{equation} \label{eq:limit 2}
    \frac{\partial}{\partial t} m_Y - \nabla_x \cdot \Bigg( D_S \nabla_x \frac{m_Y(1 + \beta m_Y)}{1 + \delta m_X + \beta m_Y} + D_H \nabla_x \frac{\delta m_X m_Y}{1 + \delta m_X + \beta m_Y} \Bigg) = - \mu m_Y + c_{HE}^{HS}\frac{\delta m_X m_Y}{1 + \delta m_X + \beta m_Y}.
\end{equation}
Apart from the mixed logistic terms appearing in equation \eqref{eq:limit 1}, the reaction--cross-diffusion system \eqref{eq:limit 1} and~\eqref{eq:limit 2} with Beddington--DeAngelis-type functional response, obtained here from a mesoscopic setting, has the same structure as the one recently derived in \cite{DesSor} using a purely macroscopic viewpoint.

\begin{remark}[The case of isotropic equilibria]
    If we consider an equilibrium state $\mathbf{M}$ with isotropic components in velocity, from the definition of the Fick coefficients \eqref{eq:diffusion coefficients} and the normalization condition \eqref{eq:property 1 of equilibria}, we easily compute the diffusion matrices as
    \begin{equation*}
        D_J = \frac{\bar{v}_J^2}{2\nu_J} I_2, \quad J \in \mathcal{S},
    \end{equation*}
    where $I_2$ denotes the two-dimensional identity matrix. Therefore, we deduce in this case that equations~\eqref{eq:limit 1} and \eqref{eq:limit 2} reduce to
    \begin{equation*}
        \begin{split}
            & \frac{\partial}{\partial t} m_X - \Delta_x \Bigg( \left(d_P \frac{1}{1 + \beta m_Y} + d_R \frac{\beta m_Y}{1 + \beta m_Y}\right) m_X \Bigg) = \alpha_P \frac{m_X}{1 + \beta m_Y} \left( 1 - \frac{c_{PP}^{P}}{\alpha_P} \frac{m_X}{1 + \beta m_Y} \right) \\[4mm]
            & \hspace{6.5cm} + \alpha_R \frac{\beta m_X m_Y}{1 + \beta m_Y} \left( 1 - \frac{c_{RR}^{R}}{\alpha_R} \frac{\beta m_X m_Y}{1 + \beta m_Y} \right) - c_{SP}^{H} \frac{m_X m_Y}{1 + \delta m_X + \beta m_Y}, \\[6mm]
            & \frac{\partial}{\partial t} m_Y - \Delta_x \Bigg( \left(d_S \frac{1 + \beta m_Y}{1 + \delta m_X + \beta m_Y} + d_H \frac{\delta m_X}{1 + \delta m_X + \beta m_Y}\right) m_Y \Bigg) = - \mu m_Y + c_{HE}^{HS}\frac{\delta m_X m_Y}{1 + \delta m_X + \beta m_Y},
        \end{split}
    \end{equation*}
    where we have denoted the diffusion coefficients $d_J = \bar{v}_J^2 / (2\nu_J)$. The left-hand sides correspond precisely to the nonlinear diffusions derived in \cite{DesSor}.
\end{remark}

%%%%%%%%%%%%%%%%%%%%%%%%%%%%%%%%%%%%%%%%%%%%%%%%%%%%
%%%%    SECTION 4: ECOLOGICAL INTERPRETATION    %%%%
%%%%%%%%%%%%%%%%%%%%%%%%%%%%%%%%%%%%%%%%%%%%%%%%%%%%

\section{Ecological interpretation of the predator-dependent prey growth}
\label{sec:ecological_interpretation}

An unexpected consequence of the kinetic derivation is the emergence of a predator-dependent carrying capacity in the prey equation. In this section, we analyze the mechanism underlying this dependence and discuss its ecological interpretation.

As already obtained in~\eqref{eq:m_P} and~\eqref{eq:m_R}, the formal fast-reaction limit between active and hidden prey yields the relations
\begin{equation}
    m_P=\frac{m_X}{1+\beta m_Y},
    \qquad
    m_R=\frac{\beta m_Xm_Y}
    {1+\beta m_Y},
\label{eq:prey_fast_equilibrium_interpretation}
\end{equation}
recalling that we have defined $\beta = c_{PS}^{RS}/c_{RE}^{P}$. Under assumptions \eqref{eq:simplified setting}, active and hidden prey have the same intrinsic growth rate $\alpha_R=\alpha_P=\alpha$ and carrying capacity $K_P=K_R=K$, while density-dependent competition acts separately within the two behavioral states. Therefore, their total contribution to prey growth is
\begin{equation*}
    \alpha m_P\left(1-\frac{m_P}{K}\right)
    +\alpha m_R\left(1-\frac{m_R}{K}\right)
    =
    \alpha m_X-\frac{\alpha}{K}\left(m_P^2+m_R^2\right)
    =
    \alpha m_X
    \left(1-\frac{m_X}{f_K(m_Y)}\right),
\label{eq:effective_logistic_growth}
\end{equation*}
where
\begin{equation}
    f_K(m_Y)
    =
K\frac{(1+\beta m_Y)^2}
    {1+(\beta m_Y)^2},
\label{eq:effective_carrying_capacity_interpretation}
\end{equation}
exploiting that $m_X=m_P+m_R$, together with relations \eqref{eq:prey_fast_equilibrium_interpretation}. Thus, the resulting prey growth is a logistic term with a predator-dependent carrying capacity $f_K(m_Y)$. Interestingly, the dependence of the prey growth term on the predator density is not introduced phenomenologically, but emerges from the redistribution of the prey population between active and refuge states.

To describe the monotonicity properties of~\eqref{eq:effective_carrying_capacity_interpretation} and its ecological interpretation,
let $q:=\beta m_Y$. A direct computation gives
\begin{equation*}
    \frac{\dd}{\dd q}
    \left(
        K\frac{(1+q)^2}{1+q^2}
    \right)
    =
    2K\frac{(1+q)(1-q)}{(1+q^2)^2}.
\end{equation*}
Consequently, $f_K$ increases for
$0<m_Y<1/\beta $, reaches its maximum when
$m_Y=1/\beta $, and decreases for~$m_Y>1/\beta $. More precisely, we have
\begin{equation*}
    f_K(0)=K,
    \qquad
    f_K\left(1/\beta \right)=2K,
    \qquad
    \lim_{m_Y\to+\infty}f_K(m_Y)=K.
    \label{eq:carrying_capacity_properties}
\end{equation*}
In particular, 
\begin{equation*}
    K\leq f_K(m_Y)\leq 2K
    \qquad\text{for every }m_Y\geq0.
\end{equation*}

The non-monotone form of $f_K$ can be interpreted in terms of intraspecific competition. In the absence of predators, almost all
prey are active, and competition is concentrated within the $P$ behavioral state, so that $m_P^2+m_R^2\simeq m_X^2$. At intermediate predator densities, the population is distributed between the active and hidden states. Since
\begin{equation*}
    m_P^2+m_R^2\leq(m_P+m_R)^2=m_X^2,
\end{equation*}
this redistribution reduces the total density-dependent competition experienced by the total prey population. The reduction is maximal when $m_P=m_R=m_X/2$, which, by
\eqref{eq:prey_fast_equilibrium_interpretation}, corresponds precisely to $\beta m_Y=1$ and hence to the maximum $f_K(m_Y)=2K$. Moreover, for very large predator densities, almost all prey are hidden, concentrating the competition within the $R$ behavioral state, and the effective carrying capacity approaches $K$ again.

This mechanism is consistent with the ecological literature on predator-induced changes in prey behavior and habitat use. Predation risk may cause prey to move from exposed foraging areas to refuges or other safer microhabitats \cite{werner1983experimental,lima1990behavioral}. Such behavioral responses may alter competitive interactions among prey and generate trait-mediated indirect effects \cite{relyea2000trait,werner2003trait,schmitz2004trophic}. In particular, experimental studies have shown that predation risk can modify, and in some cases reverse, competitive relationships through changes in habitat use and other antipredator traits \cite{relyea2000trait,jermacz2015art}. More recent field and experimental evidence also highlights the interaction between predation risk, habitat selection, and competition
\cite{creel2023habitat,peacor2025similar}.

From this perspective, $f_K(m_Y)$ should be regarded as an \emph{effective}, rather than purely environmental, carrying capacity. In fact, the presence of predators does not generate additional resources. Instead, predation risk redistributes prey between two ecological behavioral states, thereby changing how crowding and competition act on the total population. If the active and hidden states are associated with sufficiently distinct microhabitats, resources, or forms of local interference, the redistribution can alleviate competition at the population level and allow a larger total prey density to be sustained. Predator-induced movement and dispersal have indeed been shown to alter population-level density relationships and may generate non-monotonic responses \cite{cronin2024predator}.

Moreover, the positive effect of predation risk on $f_K$ distinguishes the mechanism
derived here from the classical fear effect, usually represented as a reduction in prey reproduction or resource consumption caused by antipredator behavior
\cite{peckarsky2008revisiting,wang2016modelling}. By contrast, in the present model, both active and hidden prey retain the same intrinsic growth parameters. The predator-dependent term, therefore, captures the consequences of prey redistribution, rather than a direct reproductive cost of fear. Nevertheless, the two mechanisms belong to the
broader class of non-consumptive predator effects: in both cases, predator presence changes prey population dynamics through behavioral responses rather than through direct consumption. Importantly, empirical studies indicate that such effects on prey abundance can be positive, negative, or neutral depending on ecological context, resource dynamics, and competitive interactions \cite{sheriff2020nonconsumptive,peacor2025similar}.

It is worth noting that the emergence of the effective carrying capacity in~\eqref{eq:effective_carrying_capacity_interpretation} relies on the assumption that competition acts separately within the active and hidden prey states, as specified by the reactions~{\eqref{eq:P+P to P}--\eqref{eq:R+R to R}}. This observation identifies the ecological assumption behind the new growth term: active and hidden prey must experience at least partially distinct sources of crowding or resource limitation.

Overall, the kinetic derivation provides a mechanistic justification of a predator-dependent effective carrying capacity from predator-induced prey redistribution. It complements phenomenological models with state-dependent carrying capacities \cite{liu2019dynamics} by showing how such a dependence can emerge from transitions between behavioral states and from the resulting modification of intraspecific competition.

In the following, we explore particular cases beyond assumptions \eqref{eq:simplified setting} that provide additional, interesting ecological interpretations.

%%%%    SUBSECTION: NON-REPRODUCING HIDDEN PREY    %%%%
\subsection{The case of non-reproducing hidden prey}
\label{subsec:non_reproducing_hidden_prey}

The previous derivation relies on the assumption that active and hidden prey have the same intrinsic growth rate, namely $\alpha_P=\alpha_R=\alpha$. We now consider the alternative case in which individuals entering the refuge temporarily cease to reproduce. We therefore set
\begin{equation*}
    \alpha_P=\alpha,
    \qquad
    \alpha_R=0.
\end{equation*}
The fast-reaction equilibrium is unaffected by this assumption and still gives
\begin{equation*}
    m_P=\frac{m_X}{1+q},
    \qquad
    m_R=\frac{q m_X}{1+q},
    \qquad
    q:=\beta m_Y.
\end{equation*}

It is important to distinguish whether hidden prey continue to experience density-dependent competition while in the refuge. Indeed, setting $\alpha_R=0$ removes their birth term but does not, by itself, eliminate the competition process $R+R\to R$.

\paragraph{Hidden prey subject to density-dependent competition.}

We first retain the competition term among hidden prey and assume
\begin{equation*}
    c_{PP}^{P} = \frac{\alpha}{K},
    \qquad
    c_{RR}^{R} = c_R > 0.
\end{equation*}
Notice that the parametrization $c_{RR}^{R}/\alpha_R=1/K$ used in the case $\alpha_R>0$ is no longer meaningful when~$\alpha_R=0$, hence the competition coefficient $c_R$ must be prescribed independently.

The total contribution of prey reproduction and competition is then
\begin{equation*}
    G(m_X,m_Y)
    =
    \alpha m_P-\frac{\alpha}{K}m_P^2-c_R m_R^2
    =
    \frac{\alpha m_X}{1+q}
    -
    \frac{m_X^2}{(1+q)^2}
    \left(
        \frac{\alpha}{K}+c_Rq^2
    \right).
\end{equation*}
This expression can still be written in logistic form as
\begin{equation}
    G(m_X,m_Y)
    =
    r_{\mathrm{eff}}(m_Y)m_X
    \left(
        1-\frac{m_X}{K_{\mathrm{eff}}(m_Y)}
    \right),
    \label{eq:effective_logistic_non_reproducing}
\end{equation}
where the effective growth rate is written
\begin{equation*}
    r_{\mathrm{eff}}(m_Y)
    =
    \frac{\alpha}{1+\beta m_Y}
\end{equation*}
and the effective carrying capacity becomes
\begin{equation}
    K_{\mathrm{eff}}(m_Y)
    =
    \frac{\alpha(1+q)}
    {\alpha/K+c_Rq^2}
    =
    K\frac{1+q}{1+\kappa q^2},
    \qquad \textnormal{with} \qquad
    \kappa:=\frac{Kc_R}{\alpha}.
\label{eq:effective_capacity_non_reproducing}
\end{equation}

Unlike the case $\alpha_P=\alpha_R$, predator presence now affects both the intrinsic growth rate and the effective carrying capacity. In particular, $r_{\mathrm{eff}}$ is strictly decreasing with predator density, reflecting the fact that an increasing fraction of prey is hidden and therefore does not reproduce. This provides a direct connection with classical fear-effect models, in which antipredator behavior reduces the reproductive potential of the prey population. In particular, the kinetic derivation recovers the classical reduction of prey reproduction associated with fear, but also shows that predator-induced refuge use generally modifies the density-dependent component of prey growth. 

In classical fear-effect models, predator dependence is typically introduced in the linear growth term while the quadratic competition coefficient is kept fixed. Here, behavioral redistribution also modifies the effective competition coefficient. Our derivation indicates that this separation is not generally justified when fear is mechanistically represented through transitions between reproducing active prey and non-reproducing hidden prey: the same redistribution that reduces reproduction also changes the strength of intraspecific competition.

As previously observed, the effective carrying capacity \eqref{eq:effective_capacity_non_reproducing} is non-monotone. A direct calculation gives
\begin{equation*}
    \frac{\dd}{\dd q}
    \left(
        \frac{1+q}{1+\kappa q^2}
    \right)
    =
    \frac{1-2\kappa q-\kappa q^2}
    {(1+\kappa q^2)^2}.
\end{equation*}
Consequently, it reaches its maximum at
\begin{equation*}
    q_{\max}
    =
    \sqrt{1+\frac{1}{\kappa}}-1.
\end{equation*}
At low and intermediate predator densities, the redistribution of prey between the active and hidden states may alleviate competition and increase the effective carrying capacity. At high predator densities, however, most prey accumulate in the non-reproducing refuge behavioral state while continuing to experience density-dependent competition. As a consequence, and differently from the previous (general) case, we have
\begin{equation*}
    \lim_{m_Y\to+\infty}K_{\mathrm{eff}}(m_Y)=0.
\end{equation*}

For instance, if active and hidden prey experience the same strength of density-dependent competition, that is, $c_R=\alpha/K$, then $\kappa=1$ and
\begin{equation*}
    K_{\mathrm{eff}}(m_Y)
    =
    K\frac{1+\beta m_Y}
    {1+(\beta m_Y)^2}.
\end{equation*}
The maximum is attained at $\beta m_Y=\sqrt{2}-1$, and its value is
\begin{equation*}
    K_{\mathrm{eff}}^{\max}
    =
    K\frac{1+\sqrt{2}}{2}.
\end{equation*}

\paragraph{Hidden prey without density-dependent competition.}

As a limiting scenario, one may assume that hidden prey neither reproduce nor experience appreciable density-dependent competition. This corresponds~to
\begin{equation*}
    \alpha_R=0,
    \qquad
    c_{RR}^{R}=0.
\end{equation*}
The limiting assumption $\alpha_R=0$ mathematically abstracts ecological systems where extreme refuge use severely impairs reproduction, such as the drastic reproductive suppression observed in hares under chronic predator stress~\cite{sheriff2009sensitive} or the severe metabolic constraints in hiding zooplankton~\cite{dawidowicz1992metabolic}. In this case,
\begin{equation}
    G(m_X,m_Y)
    =
    \alpha m_P
    \left(
        1-\frac{m_P}{K}
    \right)
    =
    \frac{\alpha}{1+\beta m_Y}m_X
    \left(
        1-\frac{m_X}
        {K(1+\beta m_Y)}
    \right).
\label{eq:growth_non_reproducing_noncompeting}
\end{equation}
Thus, the prey growth term can still be written in the logistic form~\eqref{eq:effective_logistic_non_reproducing} with
\begin{equation*}
    r_{\mathrm{eff}}(m_Y)
    =
    \frac{\alpha}{1+\beta m_Y},
    \qquad
    K_{\mathrm{eff}}(m_Y)
    =
    K(1+\beta m_Y).
\end{equation*}
Here the reproductive cost of refuge use produces a decreasing effective growth rate, whereas the absence of competition within the refuge produces an increasing effective carrying capacity.

Although $K_{\mathrm{eff}}$ is increasing, the net prey growth is not necessarily enhanced by predator presence. Indeed, differentiating \eqref{eq:growth_non_reproducing_noncompeting} with respect to $q=\beta m_Y$ gives
\begin{equation*}
    \frac{\partial G}{\partial q}
    =
    \frac{\alpha m_X}{(1+q)^3}
    \left(
        \frac{2m_X}{K}-(1+q)
    \right).
\end{equation*}
It follows that predator-induced refuge use decreases prey growth when
\begin{equation*}
    m_X<\frac{K}{2}(1+q),
\end{equation*}
whereas at sufficiently high prey densities the reduction in competition may temporarily compensate for the loss of reproductive individuals. This alternative scenario therefore makes explicit the balance between the reproductive cost of antipredator behavior and the competitive release generated by refuge use.

\begin{figure}[h!]
    \centering
    \includegraphics[width=\textwidth]{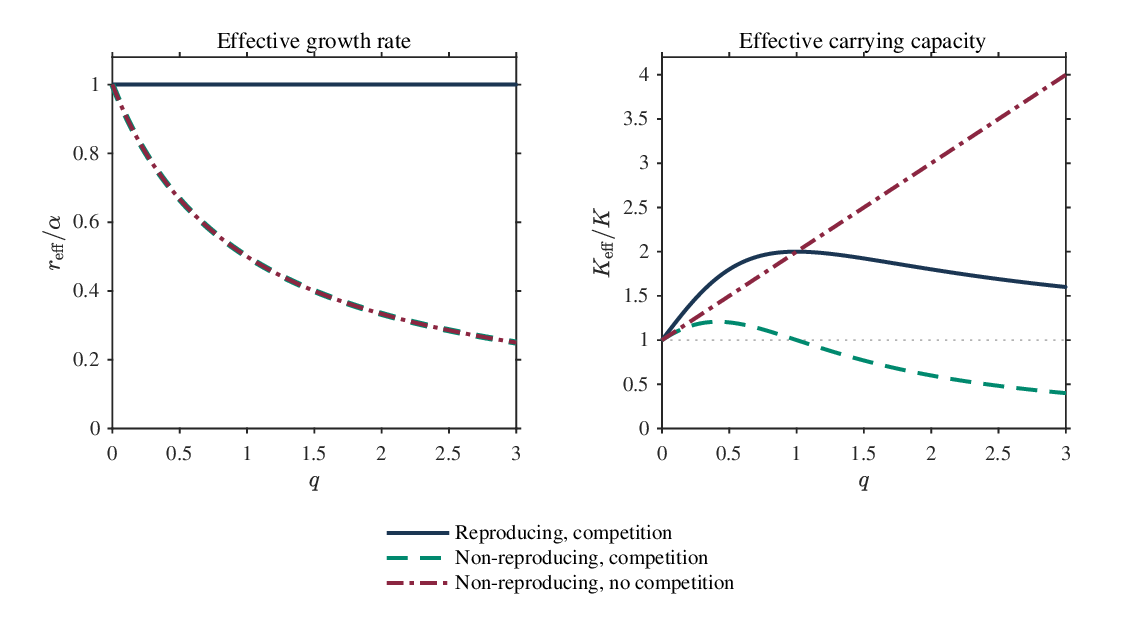}
    \caption{Dimensionless effective growth rate
    $r_{\mathrm{eff}}/\alpha$ (left panel) and effective carrying capacity
    $K_{\mathrm{eff}}/K$ (right panel) as functions of
    $q=\beta m_Y$ under the three assumptions on hidden prey reproduction and competition indicated in the legend. On the left, the curves corresponding to non-reproducing hidden prey coincide (dashed and dotted-dashed lines). On the right, the dotted horizontal line marks the predator-free carrying capacity $K$, while the curve for non-reproducing hidden prey subject to competition is plotted for $\kappa=1$, corresponding to equal competition coefficients in the active and hidden states.}
    \label{fig:effective_growth_capacity}
\end{figure}

To summarize, Figure \ref{fig:effective_growth_capacity} compares the effective growth rates and carrying capacities obtained under the three assumptions on hidden prey reproduction and competition. The left panel shows that the effective growth rate remains constant when hidden prey reproduce at the same rate as active prey. When hidden prey do not reproduce, the two curves coincide and decrease with predator density, irrespective of whether competition occurs within the hidden state.

The corresponding carrying capacities, shown on the right, exhibit qualitatively different behaviors. When both prey states reproduce and compete, the effective carrying capacity initially increases, reaches a maximum, and then approaches its predator-free value from above. When hidden prey compete but do not reproduce, the initial increase is instead followed by a decline towards zero. Finally, when hidden prey neither reproduce nor compete, the effective carrying capacity increases linearly with predator density. This increase reflects the reduced contribution of hidden prey to density-dependent losses and should therefore be interpreted together with the accompanying reduction in the effective growth rate. In this last case, the identity $r_{\mathrm{eff}}K_{\mathrm{eff}}=\alpha K$ makes this balance explicit: the carrying capacity increases while the effective growth rate tends to zero. A larger effective carrying capacity thus does not, by itself, imply faster prey population growth or a higher equilibrium abundance, which also depends on predation.

%%%%%%%%%%%%%%%%%%%%%%%%%%%%%%%%%%%%%%
%%%%    SECTION 5: CONCLUSIONS    %%%%
%%%%%%%%%%%%%%%%%%%%%%%%%%%%%%%%%%%%%%

\section{Conclusions}
\label{sec:conclusions}

In this paper, we have developed a kinetic framework for a spatial predator--prey system in which both species are structured into two behavioral states. Active and hidden prey, together with searching and handling predators, are described by distribution functions over group size, spatial position, and movement direction. The interaction mechanisms are encoded by Boltzmann-type operators, while additional linear operators account for birth and death processes. This formulation provides an explicit connection between the elementary interaction rules and the equations for the population densities, recovered as first-order moments with respect to group size. It also distinguishes the processes that redistribute individuals among behavioral states from those that change the total prey and predator populations.

The spatial description combines transport with BGK-type relaxation toward prescribed local equilibria, to describe the reorientation of species. The scaling makes explicit three distinct regimes: dominant kinetic relaxation, behavioral transitions on an intermediate time scale, and slow demographic processes. Through a moment formulation and a formal asymptotic reduction combining diffusive scaling and fast reactions, we obtain a closed macroscopic system for the total prey and predator densities. The equilibrium relations among behavioral states determine their relative contributions to both population growth and spatial movement. In particular, the dependence of these equilibrium fractions on the local population densities generates cross-diffusion when the behavioral states have different movement characteristics.

The resulting macroscopic system recovers the Beddington--DeAngelis functional response and the non-triangular diffusion structure obtained in~\cite{DesSor}. Moreover, when no state dichotomy in the prey is assumed, our approach leads to the derivation of the Holling-type II structure of \cite{ConDesSor}. The kinetic construction identifies how the transport speeds and relaxation frequencies determine the diffusivities associated with each behavioral state, and how the interaction parameters enter the limiting reaction and diffusion terms. It therefore provides a common mechanistic basis for the functional response and the cross-diffusion. Making the interaction rules and their characteristic time scales explicit also clarifies the modeling assumptions required for this macroscopic description and the origin of the additional nonlinearities produced by the reduction. The distinction between rapid behavioral transitions without prey removal and slower successful predation events clarifies the time-scale separation discussed in~\cite{iida2023cross}. In the latter process, prey removal and the associated predator-state transition occur on the same time scale.

The derivation also shows that eliminating the behavioral states modifies the effective prey-growth term. When active and hidden prey have the same intrinsic growth rate and carrying capacity, and density-dependent competition acts separately within each state, predator-induced refuge use produces a non-monotone effective carrying capacity for the total prey population. Its maximum is attained when prey are equally distributed between the two states, since this configuration minimizes aggregate competition at a given total prey density. The resulting increase in effective carrying capacity reflects a change in how competition acts on the population. Its ecological interpretation therefore depends on the assumption that the two behavioral states experience sufficiently distinct competitive interactions.

Allowing the reproductive rates to differ provides a direct connection with the fear effect. In particular, when hidden prey temporarily cease to reproduce, the fast-reaction limit yields an effective intrinsic growth rate that decreases with predator density. As far as we know, this provides the first kinetic derivation of the reproductive cost associated with the fear effect. The reduction follows from the increasing proportion of non-reproducing individuals in the refuge state. The same redistribution also changes the quadratic competition term, linking the reproductive cost of antipredator behavior to a modification of density dependence.

Within this framework, a predator-independent effective carrying capacity is not a generic consequence of refuge use. Maintaining it requires additional assumptions on how active and hidden prey contribute to density-dependent competition, so that the effective competition coefficient varies in proportion to the effective intrinsic growth rate.

The consequences depend on the competitive conditions within the refuge. If hidden prey experience negligible competition, the effective carrying capacity increases with predator density, while the effective intrinsic growth rate decreases. If hidden prey remain subject to competition, the effective carrying capacity is non-monotone and tends to zero as predator density becomes large. An increase in the effective carrying capacity alone therefore does not imply an increase in net prey growth or in equilibrium prey abundance. These outcomes show why the demographic consequences of refuge use should be assessed through the complete prey-growth term, accounting jointly for reproduction and competition.

Several questions remain open. The macroscopic reduction is formal, and a rigorous justification of the combined diffusive and fast-reaction limit would strengthen the connection between the kinetic and population-level descriptions. The limiting system also must be investigated, for instance, its equilibria and their stability, including the interaction between predator-dependent prey growth and cross-diffusion in the possible onset of spatial patterns. A further direction is to extend the kinetic construction to competing populations and investigate the emergence of the Shigesada--Kawasaki--Teramoto-type cross-diffusion~\cite{shigesada1979spatial}, which has already been justified both from microscopic~\cite{chen2021rigorous} and macroscopic~\cite{IidaMimuraNinomiya2006} models. More broadly, the present framework could be foundational for the rigorous investigation of cross-diffusion phenomena at mesoscopic scales and provides a basis for examining how behavioral transitions jointly determine effective movement, trophic interactions, and population regulation.

\section*{Acknowledgments}
All authors are members and acknowledge the support of {\it Gruppo Nazionale di Fisica Matematica} (GNFM) of {\it Istituto Nazionale di Alta Matematica} (INdAM). AB acknowledges the support from the INdAM-GNFM project CUP E5324001950001 (Multi-species non-Maxwellian Fokker-Planck models inferred from local non-equilibrium distributions), and from the European Union’s Horizon Europe research and innovation programme under the Marie Skłodowska-Curie grant agreement 101110920, project MesoCroMo (A Mesoscopic approach to Cross-diffusion Modelling in population dynamics). AB is also grateful to the University of Parma for hosting the majority of this work. CS has received funding from the Programma Giovani Ricercatori ``Rita Levi Montalcini'' 2021. MB and MG acknowledge the support of the University of Parma.\\

\section*{Disclaimer} Funded by the European Union. Views and opinions expressed are however those of the author(s) only and do not necessarily reflect those of the European Union or of the European Research Executive Agency (REA). Neither the European Union nor the granting authority can be held responsible for them. \\

%%%%%%%%%%%%%%%%%%%%%%%%%%%%
%%%%    BIBLIOGRAPHY    %%%%
%%%%%%%%%%%%%%%%%%%%%%%%%%%%

\nocite{*}
\bibliographystyle{plain}
\bibliography{Bibliography_MCD}

\end{document}